\documentclass[12pt,oneside,reqno]{amsart}
\usepackage{}
\usepackage{amssymb}

\usepackage{bbm}
\usepackage{cases}
\usepackage{amsmath}
\usepackage{graphicx}
\usepackage{mathrsfs}
\usepackage{stmaryrd}
\usepackage{amsfonts}
\usepackage{enumerate,amsmath,amssymb,amsthm}

\numberwithin{equation}{section}

\newcommand{\be}{\begin{eqnarray}}
\newcommand{\ee}{\end{eqnarray}}
\newcommand{\ce}{\begin{eqnarray*}}
\newcommand{\de}{\end{eqnarray*}}
\newtheorem{theorem}{Theorem}[section]
\newtheorem{lemma}[theorem]{Lemma}
\newtheorem{remark}[theorem]{Remark}
\newtheorem{definition}[theorem]{Definition}
\newtheorem{proposition}[theorem]{Proposition}
\newtheorem{example}[theorem]{Example}
\newtheorem{corollary}[theorem]{Corollary}

\def\e{{\mathrm{e}}}
\def\eps{\varepsilon}

\def\p{\partial}

\def\[{{\Big[}}
\def\]{{\Big]}}
\def\<{{\langle}}
\def\>{{\rangle}}
\def\({{\Big(}}
\def\){{\Big)}}

\def\bx{{\mathbf{x}}}

\def\dif{{\mathord{{\rm d}}}}

\def\min{{\mathord{{\rm min}}}}

\def\no{\nonumber}
\def\={&\!\!=\!\!&}
\def\bt{\begin{theorem}}
\def\et{\end{theorem}}
\def\bl{\begin{lemma}}
\def\el{\end{lemma}}
\def\br{\begin{remark}}
\def\er{\end{remark}}

\def\bd{\begin{definition}}
\def\ed{\end{definition}}
\def\bp{\begin{proposition}}
\def\ep{\end{proposition}}
\def\bc{\begin{corollary}}
\def\ec{\end{corollary}}
\def\bx{\begin{example}}
\def\ex{\end{example}}
\def\cA{{\mathcal A}}

\def\cI{{\mathcal I}}
\def\cJ{{\mathcal J}}

\def\cM{{\mathcal M}}
\def\cN{{\mathcal N}}

\def\cQ{{\mathcal Q}}
\def\cR{{\mathcal R}}

\def\cT{{\mathcal T}}

\def\cV{{\mathcal V}}

\def\mE{{\mathbb E}}

\def\mI{{\mathbb I}}

\def\mL{{\mathbb L}}

\def\mN{{\mathbb N}}

\def\mR{{\mathbb R}}
\def\mS{{\mathbb S}}

\def\sA{{\mathscr A}}
\def\sB{{\mathscr B}}

\def\sD{{\mathscr D}}

\def\sF{{\mathscr F}}

\def\sL{{\mathscr L}}

\def\sS{{\mathscr S}}
\def\sT{{\mathscr T}}

\def\geq{\geqslant}
\def\leq{\leqslant}

\def\div{\mathord{{\rm div}}}

\allowdisplaybreaks

\begin{document}

\title{Weak and strong well-posedness of critical and supercritical SDEs with singular coefficients}

\date{}

\author{Renming Song\ \ and \ \ Longjie Xie}

\address{Renming Song:
Department of Mathematics, University of Illinois,
Urbana, IL 61801, USA\\
Email: rsong@illinois.edu
 }

\address{Longjie Xie:
School of Mathematics and Statistics, Jiangsu Normal University,
Xuzhou, Jiangsu 221000, P.R.China\\
Email: xlj.98@whu.edu.cn
 }

\thanks{Research of R. Song is supported by the Simons Foundation (\#429343, Renming Song). L. Xie is supported by NNSF of China (No. 11701233) and NSF of Jiangsu (No. BK20170226). The Project Funded by the PAPD of Jiangsu Higher Education Institutions is also gratefully acknowledged}

\begin{abstract}
Consider the following time-dependent stable-like operator with drift
$$
\sL_t \varphi(x)=\int_{\mR^d}\big[\varphi(x+z)-\varphi(x)-z^{(\alpha)}\cdot\nabla\varphi(x)\big]\sigma(t,x,z)\nu_\alpha(\dif z)+b(t,x)\cdot\nabla \varphi(x),
$$
where $d\geq 1$, $\nu_\alpha$ is an $\alpha$-stable type L\'evy measure with $\alpha\in(0,1]$ and $z^{(\alpha)}=1_{\alpha=1}1_{|z|\leq1}z$, $\sigma$ is a real-valued Borel function on
$\mR_+\times\mR^d\times\mR^d$ and $b$ is an $\mR^d$-valued Borel function
on $\mR_+\times\mR^d$.
By using the Littlewood-Paley theory, we establish the
well-posedness for the martingale problem associated with $\sL_t$ under the sharp balance condition $\alpha+\beta\geq1$, where $\beta$ is the H\"older index of $b$ with respect to $x$.
Moreover, we also study a class of stochastic differential equations driven by Markov processes with generators of the form $\sL_t$.
We prove the pathwise uniqueness of strong solutions for such equations when the coefficients are in certain Besov spaces.
\bigskip


\end{abstract}

\maketitle \rm

\section{introduction}

Let $X=(X_t)_{t\geq 0}$ be a Feller process on $\mR^d$ and $(\sA, \sD(\sA))$ be its infinitesimal generator in $C_\infty(\mR^d)$, the space of continuous functions vanishing at infinity.
It is well known that $\sA$ satisfies the positive maximum principle, that is,
for all $\varphi\in \sD(\sA)$,
$$
\sup_{x\in\mR^d}\varphi(x)=\varphi(x_0)\geq 0\Rightarrow \sA\varphi(x_0)\leq 0.
$$
Courr\`ege's theorem then says
that, if $C_0^\infty(\mR^d)\subseteq \sD(\sA)$, then $\sA$ must be of the form
\begin{align}
\sA \varphi(x)&=\frac{1}{2}\sum_{i,j=1}^da_{ij}(x)\frac{\p^2}{\p{x_i}\p{x_j}}\varphi(x)+b(x)\cdot\nabla \varphi(x)\no\\
&\quad+\int_{\mR^d}\big[\varphi(x+z)-\varphi(x)-1_{\{|z|\leq 1\}}z\cdot\nabla \varphi(x)\big]\mu(x,\dif z),\quad\forall \varphi\in C^{\infty}_0(\mR^d),\label{ge}
\end{align}
where for each $x\in\mR^d$, $a(x)=(a_{i,j}(x))$ is a
non-negative definite real symmetric $d\times d$-matrix,
$b(x)$ is an $\mR^d$-valued function,
and $\mu(x,\cdot)$ is a measure on $\mR^d\setminus\{0\}$ such that
$$
\int_{\mR^d\setminus\{0\}}\big(1\wedge|z|^2\big)\mu(x,\dif z)<\infty.
$$
The triplet $(a, b, \mu)$ is uniquely determined by $\sA$. In this way, all information about the process $X$ are contained in the operator or the triplet $(a, b, \mu)$.
However, the opposite question that whether a given operator $\sA$ of the form (\ref{ge}) (or a given triplet $(a, b, \mu)$)
actually generates a unique Markov process in $R^d$ with $\sA$ as its infinitesimal generator
is pretty difficult. Many people have studied this question.

\vspace{2mm}
Let us consider the following more general time-dependent L\'evy type operator: for every $\varphi\in C^{\infty}_0(\mR^d)$,
\begin{align}
\sA_t^{a, b, \mu} \varphi(x)&
:=\frac{1}{2}\sum_{i,j=1}^da_{ij}(t,x)\frac{\p^2}{\p{x_i}\p{x_j}}\varphi(x)+b(t,x)\cdot\nabla \varphi(x)\no\\
&\quad+\int_{\mR^d}\big[\varphi(x+z)-\varphi(x)-1_{\{|z|\leq 1\}}z\cdot\nabla \varphi(x)\big]\mu(t,x,\dif z).\label{ge2}
\end{align}
There are many ways to specify
a Markov process corresponding to $\sA_t^{a, b, \mu}$,
and the common ones are: Stroock-Varadhan's martingale problem which determines the distribution of the process,  and It\^o's the stochastic differential equation (SDE) which
gives a process as a strong solution, or more generally a weak solution.
From a probability point of view, the three parts on the right hand side of (\ref{ge2}) are, respectively, the diffusion term, the drift term and the jump term.
The second order differential operator $\sA_t^{a, b, 0}$ and the corresponding diffusion process
have been intensively studied both in probability and partial differential equations.
Under mild conditions on the coefficients $a$ and $b$,
it can be shown, using the martingale problem method (see \cite{S-V}), that
there exists a unique diffusion process $X$ having $\sA_t^{a, b, 0}$ as its infinitesimal generator.  This diffusion process can also be constructed as a weak solution to the following It\^o's stochastic differential equation:
\begin{align}
\dif X_t=\sigma(t,X_t)\dif W_t+b(t,X_t)\dif t,\quad X_0=x\in\mR^d,\label{sde0}
\end{align}
where $a=\sigma\sigma^*$ and $W_t$ is a standard Brownian motion.
We also mention that,  without the Lipschitz condition, Veretennikov \cite{Ve}
proved that  \eqref{sde0} has a unique strong solution
when $\sigma=\mI$ and $b\in L^\infty(\mR_+\times\mR^d)$.
A further extension was obtained by  Krylov and R\"ockner \cite{Kr-Ro} where the pathwise uniqueness for (\ref{sde0}) was shown under the condition that
$$
b\in L^q\big(\mR_+;L^p(\mR^d)\big)\quad\text{with}\quad d/p+2/q<1.
$$
Note that under these conditions, the corresponding deterministic ordinary differential system (i.e., $\sigma=0$) is far from being well-posed. This is usually called, following the terminology of Flandoli \cite{Fla}, a regularization by noise phenomenon.

\vspace{2mm}
Nowadays, much attention has been paid to non-local operators and their corresponding discontinuous processes, due to their importance both in theory and in applications. Consider the following time-independent non-local operator
$$
\sA_\mu\varphi(x):=\int_{\mR^d}\big[\varphi(x+z)-\varphi(x)-1_{\{|z|\leq 1\}}z\cdot\nabla \varphi(x)\big]\mu(x,\dif z).
$$
The simplest circumstance would be $\mu(x,\dif z)\equiv c_{d,\alpha}|z|^{-d-\alpha}\dif z$ with $\alpha\in(0,2)$ and $c_{d,\alpha}>0$ being a constant.
In this case, $\sA_\mu$ is the fractional Laplacian  $\Delta^{\alpha/2}$ which is the generator of an isotropic $\alpha$-stable process.
The next natural generalization would be
$\alpha$-stable like operators, that is, the case where the jump measure $\mu$ is state-dependent, absolutely continuous with respect to the Lebesgue measure with a jump intensity comparable to that of an  $\alpha$-stable process. In the literature, there are two different meanings to the term ``$\alpha$-stable like operators'': one refers to (see \cite{Ba})
$$
\sA_{\alpha(x)}\varphi(x):=\int_{\mR^d}\big[\varphi(x+z)-\varphi(x)-1_{\{|z|\leq 1\}}z\cdot\nabla \varphi(x)\big]
\frac{c_{d,\alpha(x)}}{|z|^{d+\alpha(x)}}\dif z,
$$
and the other refers to (see \cite{CK})
\begin{align}
\sA_{\kappa}\varphi(x):=\int_{\mR^d}\big[\varphi(x+z)-\varphi(x)-1_{\{\alpha\geq1\}}1_{\{|z|\leq 1\}}z\cdot\nabla \varphi(x)\big]\frac{\kappa(x,z)}{|z|^{d+\alpha}}\dif z. \label{ak}
\end{align}
The martingale problem for $\sA_{\alpha(x)}$ and $\sA_\kappa$ have been studied in \cite{Ba} and \cite{B-T}, respectively, see also \cite{BBCK,B-K-K,Ch-Zh,C-Z} for
related results and references therein.
We also mention that in \cite{B-C}, the weak uniqueness for a system of SDEs driven by a cylindrical $\alpha$-stable process was studied, that is, for $i=1,\cdots,d$,
$$
\dif X_t^i=\sum_{j=1}^da_{i,j}(X_{t-})\dif L_t^j,\quad X_0=x\in\mR^d,
$$
where $L_t^j$ are independent one-dimensional symmetric $\alpha$-stable processes. Note that in this case, the L\'evy measure $\nu$ of the process $\hat L_t:=(L_t^1,\cdots,L_t^d)$ is given by
\begin{align}
\nu(\dif z)=\sum_{j=1}^d1_{\{|z_1|=0,\cdots,|z_j|\neq0,\cdots,|z_d|=0\}}|z_j|^{-d-\alpha}\dif z_j,\quad\forall z=(z_1,\cdots,z_d)\in\mR^d,   \label{cyl}
\end{align}
which is not absolutely continuous with respect to the Lebesgue measure.

\vspace{2mm}
When the non-local operator $\sA_\mu$ is perturbed  by a gradient operator, the situation is much more complicated.
Let us consider the following fractional Laplacian  with gradient perturbation: for $\alpha\in(0,2)$,
\begin{align*}
\sA^b_\alpha \varphi(x):=\Delta^{\alpha/2}\varphi(x)+b(t,x)\cdot\nabla \varphi(x),
\end{align*}
which corresponds to the following SDE driven by an isotropic $\alpha$-stable process
\begin{align}
\dif X_t=\dif L_t+b(t,X_t)\dif t,\quad X_0=x\in\mR^d.\label{sde00}
\end{align}
Even in this simplest case, the study of the operator $\sA^b_\alpha$ and SDE (\ref{sde00}) is much more delicate than that of $\sA^{a, b, 0}$ and SDE (\ref{sde0}) due to the non-local nature of the generator and the discontinuity of the process. In fact,
in the case $\alpha\in(1,2)$,
the non-local operator $\Delta^{\alpha/2}$ is the dominant term and $b\cdot\nabla$ can be seen as a lower order perturbation of
$\Delta^{\alpha/2}$. Thus the case $\alpha\in(1,2)$ is called the subcritical case.
The critical case corresponds to $\alpha=1$
since $\Delta^{1/2}$ and $\nabla$ are of the same order.
The case $\alpha\in(0,1)$ is called the supercritical case since in this case the gradient term $b\cdot\nabla$ is of higher order than $\Delta^{\alpha/2}$.

Up to now, the martingale problem for $\sA^b_\alpha$ as well as the weak and strong well-posedness for SDE (\ref{sde00}) in the subcritical case have also been intensively studied. We refer to \cite{B} for an overview and a rich reference list, see also \cite{K-S}.
Consider the following more general $\alpha$-stable like operator perturbed by gradient:
$$
\sA_\kappa^b\varphi(x)=\sA_\kappa\varphi(x)+b(t,x)\cdot\nabla\varphi(x),
$$
where $\sA_\kappa$ is defined by (\ref{ak}).
The martingale problem for $\sA_\kappa^b$
with $\alpha\in(1,2)$ and bounded drift was studied in \cite{M-P}, see also \cite{A_K,H,S-W}  for related results using the theory of pseudo-differential operators.
We also mention that
Priola \cite{Pri}  proved that SDE (\ref{sde00}) admits a
pathwise unique strong solution when $\alpha\geq1$ and
$b(t,x)=b(x)\in C_b^\beta(\mR^d)$ with $\beta>1-\alpha/2$.
Later, Zhang \cite{Zh00} obtained the strong well-posedness for SDE (\ref{sde00}) when $\alpha>1$ and
$$
b\in L^\infty\big(\mR_+;L^\infty(\mR^d)\cap H^\beta_p(\mR^d)\big)\quad\text{with}\quad\beta>1-\alpha/2\quad\text{and}\quad p>2d/\alpha.
$$
See also \cite{X-X,X-Z} for generalizations to the multiplicative noise cases.

\vspace{2mm}
The critical and supercritical cases
are more difficult and there are considerably less results in the literature in these two cases.
In fact, when $d=1$ and $L_t$ is an isotropic $\alpha$-stable process with $\alpha<1$, a counterexample was given by Tanaka, Tsuchiya and Watanabe \cite{Ta-Ts-Wa} which showed that the weak uniqueness for SDE (\ref{sde0}) fails even if $b$ is bounded, time-independent and $\beta$-H\"older
continuous with
$$
\alpha+\beta<1,
$$
see also \cite{B-B-C}.
On the other hand, the weak uniqueness for the one-dimensional SDE (\ref{sde0}) was obtained in  \cite{Ta-Ts-Wa} with non-decreasing $\beta$-H\"older continuous drift under the additional condition
$$
\alpha+\beta>1.
$$
In view of this, we shall call $\alpha+\beta=1$ the balance condition below.
Very recently, Chen, Song and Zhang \cite{Ch-So-Zh} studied the strong well-posedness of SDE (\ref{sdee})  when $0<\alpha\leq 1$ and
$$
b\in L^\infty\big(\mR_+;C_b^\beta(\mR^d)\big)\quad\text{with}\quad \beta>1-\alpha/2.
$$
This was generalized to the multiplicative noise case with Lipschitz jump diffusion coefficient in \cite{C-Z-Z}.

\vspace{2mm}

For any Borel function $\sigma(t,x,z): \mR_+\times\mR^d\times\mR^d\to\mR$ and L\'evy measure $\nu$, we define
\begin{align}
\sL^\sigma_{\nu}\varphi(x):=\int_{\mR^d}\big[\varphi(x+z)-\varphi(x)
-1_{\{\alpha\geq 1\}}1_{\{|z|\leq 1\}}z\cdot \nabla \varphi(x)\big]\sigma(t,x,z)\nu(\dif z).\label{ope}
\end{align}
In this paper, we consider the following time-dependent non-local and non-symmetric L\'evy type operator:
\begin{align}
\sL_t\varphi(x):=\sL_\nu^\sigma\varphi(x)+b(t,x)\cdot\nabla\varphi(x),  \quad\forall \varphi\in C_0^\infty(\mR^d),\label{oper}
\end{align}
where $b(t,x):\mR_+\times\mR^d\to \mR^d$ is a Borel  function.
We will assume that $\nu$ is of $\alpha$-stable type,
which may be singular with respect to the Lebesgue measure, and contains the cylindrical case (\ref{cyl}) as an example.
Moreover, we will concentrate on  the case $\alpha\in(0,1]$, that is, the supercritical and critical cases.

The first aim of this paper is to prove the uniqueness for martingale problem associated with the operator $\sL_t$ under the sharp balance condition $\alpha+\beta\geq1$, where $\beta$ is the H\"older index of the drift $b$ with respect to $x$, see Theorems \ref{main1} and  \ref{main2}.
The unique martingale solution for $\sL_t$ is known to be equivalent to the weak solution of a class of stochastic differential equations. The second aim of this paper is to study the strong well-posedness of such equations under weak assumptions on the coefficients. We point out that due to the
state-dependent jump intensity kernel $\sigma(t,x,z)$, the usual SDE driven by a L\'evy process is not suitable to characterize $\sL_t$ (see \cite{Xie}). To specify the SDE we are going to study, let
$\cN(\dif z,\dif r,\dif t)$ be a Poisson random measure on $\mR^d\times[0,\infty)\times[0,\infty)$ with intensity measure $\nu(\dif z)\dif r\dif t$, and
$\widetilde \cN(\dif z,\dif r,\dif t):=\cN(\dif z,\dif r,\dif t)-\nu(\dif z)\dif r\dif t$
be the corresponding compensated Poisson random measure.
Then the Markov process $X_t$ corresponding to $\sL_t$ should satisfy the following SDE:
\begin{align}
\dif X_t&=\int_0^{\infty}\!\!\!\!\int_{|z|\leq 1}\!\!1_{[0,\sigma(s,X_{s-},z)]}(r)z\widetilde \cN(\dif z,\dif r,\dif t)\no\\
&\quad+\int_0^{\infty}\!\!\!\!\int_{|z|> 1}\!\!1_{[0,\sigma(s,X_{s-},z)]}(r)z \cN(\dif z,\dif r,\dif t)+b(t,X_t)\dif t, \ \ \ \ X_0=x\in\mR^d. \label{sdee}
\end{align}
In fact, noticing that for any  function $f$ on $\mR^d$ and any $r>0$,
\begin{align*}
f\big(x+1_{[0,\sigma(t,x,z)]}(r)z\big)-f(x)=1_{[0,\sigma(t,x,z)]}(r)\big[f(x+z)-f(x)\big],
\end{align*}
an application of It\^o's formula shows that the infinitesimal generator of the solution to SDE (\ref{sdee}) is given exactly by (\ref{oper}). Note also that the driving noise is a Markov process which is not necessarily of L\'evy type (see \cite{Kurz2}).
Under the conditions that $b$ is bounded and globally Lipschitz, $\sigma$ is bounded with
\begin{align*}
\int_{\mR^d}|\sigma(x,z)-\sigma(y,z)|\cdot|z|\nu(\dif z)\leq C_1|x-y|, \ \ \ \ \forall \ \ x, y \in \mR^d,
\end{align*}
and some other assumptions, Kurtz \cite{Kurz2} proved the existence and uniqueness of strong solution to SDE \eqref{sdee}, see also \cite{B-G-R,Im-Wi} and references therein for related results and applications to Boltzmann equations.
We will prove the strong well-posedness of  SDE \eqref{sdee} with coefficients in certain Besov spaces (see Definition \ref{bes}), see Theorem \ref{main3}.
Our results show that, even in the critical and supercritical case, the phenomenon of regularization by noise happens.

\vspace{2mm}
Now we give a brief description of the strategy of our proofs.
Both the weak well-posedness and the strong well-posedness
for SDE (\ref{sdee}) rely on  the following supercritical drift-diffusion equation: for every $\lambda\geq0$,
\begin{align}
\p_tu=\sL_\nu^\sigma u+b\cdot\nabla u-\lambda u+f.    \label{eqq}
\end{align}
Such an equation is
of independent interest in itself  since it is closely related to quasi-geostrophic equation, hydrodynamic transport equation and 3-D Navier-Stokes equation, see \cite{C-V,CMZ,C-W,M-M,Si} for the study of (\ref{eqq}) when $\sL_\nu^\sigma\equiv\Delta^{\alpha/2}$ with $\alpha\in(0,1]$.
We will use the Littlewood-Paley theory to
obtain the optimal regularity of (\ref{eqq}) both in Besov spaces and in H\"older spaces, see Theorems \ref{drift3}, \ref{drift} and  \ref{th}.
Then, the uniqueness of the martingale problem for $\sL_t$ will
follow by an application of It\^o's formula. However, as we will see,  under our assumptions, the solution $u$ will not be regular enough to apply the It\^o's formula. To overcome this difficulty, we use a commutator estimate in the supercritical case $\alpha\in(0,1)$ and Krylov's estimate in the critical case $\alpha=1$.
As for the strong uniqueness for SDE (\ref{sdee}) with singular coefficients, we will adopt Zvonkin's  argument to transform SDE (\ref{sdee}) into a new one with better coefficients. Krylov's estimate again will play a key role, see Lemma \ref{kkr}. We also point out that, unlike the classical SDEs driven by multiplicative L\'evy noise considered in \cite{Ch-So-Zh,C-Z,C-Z-Z,Kr-Ro,Pri,X-Z,Zh00},  the usual $L_2$-difference-estimate is not applicable for studying SDE (\ref{sdee}) due to the extra term $1_{[0,\sigma(X_{s-},z)]}(r)$. Instead, we willll use a mixture of $L_1$ and $L_2$ estimates, see \cite{Kurz2}. Due to the irregularity of $b$ and $\sigma$, it is much more complicated than in \cite{Kurz2} to apply this technique.

\vspace{2mm}

The organization of this paper is as follows: In Section 2, we state our main results, including the weak and strong well-posedness for SDE \eqref{sdee}.
We present in Section 3 some preliminaries of
the Littlewood-Paley theory. Section 4 is devoted to study the regularities of
the parabolic integral-differential equation (\ref{eqq}) both in Besov spaces and in H\"older spaces. Finally, the main results are proved in Section 5.

We conclude this section by spelling out some conventions that will be used throughout this paper. For $a, b\in \mR$, $a\vee b:=\max\{a, b\}$ and $a\wedge b:=\min\{a, b\}$.
For any $R>0$, $B_R:=\{x\in\mR^d:|x|\leq R\}$.
The letter $C$ with or without
subscripts will denote an unimportant constant, whose value may change in different places, and whose dependence on parameters can be traced from calculations. We write $f\preceq g$ to
mean that $f\leq Cg$ for some $C>0$. We will use $:=$ to denote a definition, and we assume that all the functions
considered in this paper are Borel.

\section{Statement of main results}

We first specify the conditions that we will impose on the coefficients $\sigma$, $b$ and the L\'evy measure $\nu$ of the operator $\sL_t$ given by \eqref{oper}.
For $\alpha\in(0,2)$, we denote by $\mL^\alpha$ the space of all non-degenerate
symmetric $\alpha$-stable L\'evy measures $\nu^\alpha$, that is,
\begin{align*}
\nu^\alpha(A)=\int_0^\infty\!\left(\int_{\mS^{d-1}}\!1_A(r\theta)\Sigma(\dif \theta)\right)\frac{\dif r}{r^{1+\alpha}},\quad \forall A\in \sB(\mR^d),
\end{align*}
where $\Sigma$ is a finite symmetric  measure on the unit sphere $\mS^{d-1}$ in $\mR^d$ satisfying
\begin{align}\label{nu}
\int_{\mS^{d-1}}\!|\theta_0\cdot\theta|^\alpha\Sigma(\dif \theta)>0,\quad\forall\theta_0\in\mS^{d-1}.
\end{align}
We will use $\psi_{\nu^\alpha}$ to denote the symbol of the purely discontinuous L\'evy process
with L\'evy measure $\nu^\alpha$ which is given by
$$
\psi_{\nu^\alpha}(\xi)=\int_{\mR^d}\big(\e^{i\xi\cdot y}-1
-1_{\{\alpha\geq 1\}}1_{\{|y|\leq 1\}}i\xi\cdot y\big)\nu^\alpha(\dif y).
$$
Then, it is known that (see \cite{Pri}) condition \eqref{nu} is also equivalent to the condition that for some constant $c_0>0$,
\begin{align}\label{cy}
\psi_{\nu^\alpha}(\xi)\leq -c_0|\xi|^\alpha,\quad\forall\xi\in\mR^d.
\end{align}
Throughout this paper, we make the following assumption on the L\'evy measure $\nu$:

\vspace{1mm}
\begin{enumerate}
\item [{\bf (H$^\nu$)}]
$\nu(B^c_1)<\infty$ and there exist $\nu_1, \nu_2\in\mL^\alpha$ such that
\begin{align*}
\nu_1(A)\leq \nu(A)\leq \nu_2(A),\quad \forall A\subseteq B_1.
\end{align*}
\end{enumerate}

\br
i) The mixed-stable case $\nu(\dif z)=
(c_{d,\alpha}|z|^{-d-\alpha}+c_{d,\beta}|z|^{-d-\beta})\dif z$
with $\beta<\alpha$
as well as the truncated $\alpha$-stable-like case $\nu(\dif z)=c_{d,\alpha}1_{\{|z|\leq1\}}\kappa(z)|z|^{-d-\alpha}\dif z$ with $0\leq\kappa_0\leq\kappa(z)\leq\kappa_1$ satisfy {\bf (H$^\nu$)}.

ii) Note the the L\'evy measure $\nu$ can be singular with respect to the Lebesgue measure. In particular, for $d\geq2$, it is easy to see that the cylindrical $\alpha$-stable case (\ref{cyl}) satisfies (\ref{cy}), and hence {\bf (H$^\nu$)}.
\er

Here are some possible assumptions that we will make on the coefficients $\sigma$ and $b$:

\vspace{1mm}
\begin{enumerate}
\item [{\bf (H$^\sigma_1$)}] There exist constants $\kappa_0,\kappa_1>0$, $\kappa_2\geq1$ and $\vartheta\in(0,1]$ such that
\begin{align*}
\kappa_0\leq \sigma(t,x,z)\leq \kappa_1,
\quad \forall (t,x,z)\in\mR_+\times\mR^d\times\mR^d
\end{align*}
and for all $(t, z)\in\mR_+\times\mR^d$,
\begin{align*}
|\sigma(t,x,z)-\sigma(t,y,z)|\leq \kappa_2|x-y|^\vartheta,
\quad \forall x, y\in \mR^d \mbox{ with } |x-y|\leq  1.
\end{align*}
\item [{\bf (H$^\sigma_2$)}] There exists a function $\varrho\in B^0_{q,\infty}(\mR^d)$ with $q>d/\alpha$ such that for every $t>0$ and almost all $x,y\in \mR^d$,
\begin{align}
\int_{\mR^d}|\sigma(t,x,z)-\sigma(t,y,z)|(|z|\wedge1)\nu(\dif z)\leq |x-y|\Big(\varrho(x)+\varrho(y)\Big).\label{a1}
\end{align}

\item [{\bf (H$^b_1$)}] $b\in L^\infty\big(\mR_+;B^{\beta}_{p,\infty}(\mR^d)\big)$
with
$$
\beta>1-\alpha\quad\text{and}\quad d/(\alpha+\beta-1)<p\leq\infty.
$$

\item [{\bf (H$^b_2$)}] For $\alpha\in(0,1)$, $b\in L^\infty\big(\mR_+;C^{1-\alpha}_{b}(\mR^d)\big)$
with $\|b\|_{L^\infty C^{1-\alpha}_{b}}$ small; for $\alpha=1$, $b\in L^\infty(\mR_+\times\mR^d)$
with  $\|b\|_{L^\infty}$ small.

\item [{\bf (H$^b_3$)}] $b\in L^\infty\big(\mR_+;B^{\beta}_{p,\infty}(\mR^d)\big)$
with
$$
\beta>1-\alpha/2\quad\text{and}\quad 2d/\alpha<p\leq\infty.
$$
\end{enumerate}

Our main results concerning the well-posedness of the martingale problem for $\sL_t$ and the weak uniqueness of SDE (\ref{sdee}) are as follows.

\bt\label{main1}
Let $0<\alpha\leq 1$. Assume that
{\bf (H$^\sigma_1$)} and {\bf (H$^b_1$)} hold.
Then SDE (\ref{sdee}) has a unique weak solution for every $x\in \mR^d$.
Equivalently, the martingale problem for $\sL_t$ has a unique solution.
\et

\br
Note that by (\ref{bh}) below, we have $B^\beta_{\infty,\infty}(\mR^d)=C_b^\beta(\mR^d)$. Thus, Theorem \ref{main1} not only generalizes the main result of \cite{Ta-Ts-Wa} to the multidimensional case, but also to the case with more general variable coefficients.
\er

\bt\label{main2}
Let $0<\alpha\leq1$. Assume that
{\bf (H$^\sigma_1$)} and {\bf (H$^b_2$)} hold.
Then SDE (\ref{sdee}) admits a unique weak solution for every $x\in \mR^d$.
Equivalently, the martingale problem for $\sL_t$ has a unique solution.
\et

\br
In the critical case $\alpha=1$, a similar result under the balance condition was proved by Tsuchiya \cite{Ts} in dimension 1. That is, the uniqueness of the martingale problem for $\Delta^{1/2}+b(x)\cdot\nabla$ in $\mR^1$ was proved when $b$ is bounded with $\|b\|_\infty$ small. Thus, we not only generalize
the main result of \cite{Ts} to higher dimensions in the critical case $\alpha=1$
but also to the supercritical case $\alpha\in(0,1)$. In view of \cite{B-B-C,Ta-Ts-Wa}, our result are almost sharp.
However, for general bounded and $(1-\alpha)$-H\"older drift $b(t,x)$, the problem is still open.
We also point out that if the constant $\kappa_0$ in {\bf (H$^\sigma_1$)} is large enough, then the smallness condition on $b$ can be dropped, see Remark \ref{ppp}.
\er

We also study the strong well-posedness of SDE (\ref{sdee}) with irregular coefficients and prove the following result.

\bt\label{main3}
Let $0<\alpha\leq 1$. Assume that
{\bf (H$^\sigma_1$)},
{\bf (H$^\sigma_2$)} and {\bf (H$^b_3$)} hold.
Then for each $x\in\mR^d$, there exists a unique strong solution $X_t(x)$ to SDE (\ref{sdee}).
\et

\br
i) Compared with \cite{Ch-So-Zh,C-Z-Z}, our assumptions are weaker both in the
drift coefficient $b$  and in the coefficient $\sigma$.
Moreover, SDE (\ref{sdee}) is more general than SDE (\ref{sde00}).

ii) Note  that it holds that $L^q(\mR^d)\subset B^0_{q,\infty}(\mR^d)$. Thus, if $\sigma(t,x,z)\equiv g(t,x)$ with $g\in L^\infty\big(\mR_+;W^{1,p}(\mR^d)\big)$ and $p>d/\alpha$, where $W^{1,p}(\mR^d)$ is the usual first order Sobolev space, then {\bf (H$^\sigma_2$)} is satisfied. Moreover,  by the Sobolev embedding (\ref{emb2}) below, the assumption on the drift $b$ in {\bf (H$^b_3$)}  is in fact the same if we replace $B^{\beta}_{p,\infty}(\mR^d)$ with $W^{\beta,p}(\mR^d)$ or $H^{\beta,p}(\mR^d)$.
\er

\section{Preliminaries}

We first recall some preliminaries of the Littlewood-Paley theory. For more details, see e.g., \cite{B-C-D}.
Let $\sS(\mR^d)$ be the Schwartz space of all rapidly decreasing functions
and $\sS'(\mR^d)$ be its dual space which consists of all tempered distributions.
Given  $f\in\sS'(\mR^d)$, we denote by $\sF f=\hat f$ (resp. $\sF^{-1}f=\check{f}$) the Fourier transform (resp. the Fourier inverse transform) of $f$.
The following definition is well known.

\bd
A dyadic partition of unity $(\chi, \rho)$ consists of two smooth functions $\chi, \rho\in C^\infty(\mR^d)$ taking values in $[0, 1]$ such that ${\rm supp}\chi\subseteq B_2$, ${\rm supp}\rho\subseteq B_{2}\backslash B_{1/2}$,
and satisfying that
\begin{align*}
\chi(x)+\sum_{j\geq0}\rho(2^{-j}x)=1,&\quad \forall x\in\mR^d.
\end{align*}
We shall write $\rho_{-1}(x):=\chi(x)$ and $\rho_j(x):=\rho(2^{-j}x)$ for all $j\geq0$.
\ed

From now on, we fix such a dyadic partition of unity $(\chi,\rho)$ and define the Littlewood-Paley operators as follows: for each $f\in\sS'(\mR^d)$,
\begin{align*}
\Lambda_jf:=0\quad{{\rm if}}\quad j\leq-2,\quad{\rm and}\quad\Lambda_jf:=\rho_j(D)f:=\sF^{-1}(\rho_j\sF f)\quad{\rm if}\quad j\geq-1.
\end{align*}
Informally, $\Lambda_j$ is a frequency projection to the annulus $\{|\xi|\approx2^j\}$. We also introduce the low-frequency cut-off operator
$$
S_jf:=\sum_{i\leq j-1}\Lambda_if,
$$
as well as $h_j:=\sF^{-1} \rho_j$ so that
\begin{align}
h_j\ast f=\sF^{-1}(\rho_j\sF f)=\Lambda_jf,  \label{con}
\end{align}
where $\ast$ denotes the usual convolution.
Then, it is known that
\begin{align}
f&=\sum_{j\geq -1}\Lambda_jf=\lim_{j\rightarrow\infty}S_jf,\quad \forall f\in\sS'(\mR^d),    \label{id}
\end{align}
where the limit is taken in the space $\sS'(\mR^d)$.
Notice that with our choice of the dyadic partition of unity, it is easy to verify that
\begin{equation} \label{00}
\begin{aligned}
\Lambda_j\Lambda_kf=0\quad&\text{if}\quad|j-k|\geq2,\\
\Lambda_j(S_{k-1}f\Lambda_kg)=0\quad&\text{if}\quad|j-k|\geq5.
\end{aligned}
\end{equation}
For  $f, g\in \sS'(\mR^d)$, the paraproduct of $g$ by $f$ is defined by
$$
\cT_fg=\sum_{i}S_{i-1}f\Lambda_ig,
$$
and the remainder of $f$ and $g$ is defined by
$$
\cR(f,g):=\sum_{i}\sum_{|j|\leq1}\Lambda_if\Lambda_{i-j}g.
$$
Then, we have the following  Bony decomposition of product:
\begin{align}
fg=\cT_fg+\cT_gf+\cR(f,g).    \label{bon}
\end{align}

Now we recall the definition of Besov spaces.
\bd\label{bes}
For $\beta\in\mR$ and $p,q\in[1,\infty]$, the nonhomogeneous Besov space $B^\beta_{p,q}(\mR^d)$ is defined as the set of all $f\in\sS'(\mR^d)$ such that
$$
\|f\|_{B^\beta_{p,q}}:=1_{\{q<\infty\}}\left(\sum_{j\geq-1}\big(2^{\beta j}\|\Lambda_jf\|_p\big)^q\right)^{1/q}+1_{\{q=\infty\}}\left(\sup_{j\geq-1}2^{\beta j}\|\Lambda_jf\|_p\right)<\infty,
$$
here and below, $\|\cdot\|_p$ denotes the usual $L^p$-norm in $\mR^d$.
\ed

Let us list some elementary properties of Besov spaces which will be used later.
It is known that for $\beta\in(0,\infty)\setminus\mN$, we have
\begin{align}
B^\beta_{\infty,\infty}(\mR^d)=C^\beta_b(\mR^d),    \label{bh}
\end{align}
where  $C^{\beta}_b(\mR^d)$ is the usual H\"older space consisting of functions whose
partial derivatives of order $[\beta]$ are $(\beta-[\beta])$-H\"older continuous. We emphasize that in the case $\beta\in\mN$, the space $B^{\beta}_{\infty,\infty}(\mR^d)$ is strictly larger than $C^\beta_b(\mR^d)$. For $\beta>0$ which is not an integer and $1<p<\infty$, $1\leq q\leq\infty$, it holds that
\begin{align}
B^{\beta+\frac{d}{p}}_{p,q}(\mR^d)\hookrightarrow C_b^\beta(\mR^d).    \label{embb}
\end{align}
Note that for $p\geq2$,
\begin{align}
B^0_{p,1}(\mR^d)\subseteq L^p(\mR^d)\subseteq B^0_{p,\infty}(\mR^d).  \label{pp}
\end{align}
We also have the following embedding relations between Besov spaces: for any $\beta\in\mR$ and $p,q\in[1,\infty]$, it holds that
\begin{align}
B^\beta_{p,q}(\mR^d)\subseteq B^\beta_{p,\infty}(\mR^d),   \label{emb0}
\end{align}
and for any $\beta_1,\beta_2\in\mR$ and $p_1,q_1\in[1,\infty]$ with
$$
p\leq p_1,\,\,\,q\leq q_1,\,\,\,\beta<\beta_2\quad\text{and}\quad \beta-d/p=\beta_1-d/p_1,
$$
it holds that
\begin{align}
B^{\beta_2}_{p,\infty}(\mR^d)\subseteq B^\beta_{p,q}(\mR^d)\subseteq B^{\beta_1}_{p_1,q_1}(\mR^d).    \label{emb}
\end{align}
Below, for $\theta\in[0,1]$ and two Banach spaces $A$, $B$, we use $[A, B]_{\theta}$ to denote the complex interpolation space between $A$ and $B$. It is well known that there is a constant $c_{\theta}>0$ such that
\begin{align}
\|y\|_{[A, B]_{\theta}}\leq c_{\theta}\|y\|_{A}^{1-\theta}\|y\|_{B}^{\theta},\quad\forall y\in A\cap B.  \label{ineq}
\end{align}
For the interpolation between Besov spaces, we have that for $\beta_0,\beta_1\in\mR$, $p>1$ and any $\theta\in(0,1)$,
\begin{align}
[B^{\beta_0}_{p,\infty}(\mR^d), B^{\beta_1}_{p,\infty}(\mR^d)]_\theta=B^{\beta}_{p,\infty}(\mR^d), \label{int}
\end{align}
where $\beta=(1-\theta)\beta_0+\theta\beta_1$. The above facts are standard and can be found in \cite{B-C-D} or \cite{Tri}.

For $0<\beta\leq2$, we also recall the Bessel potential space $H^{\beta,p}(\mR^d)$ which is given by
$$
H^{\beta,p}(\mR^d):=
\left\{f\in L^p(\mR^d): \Delta^{\beta/2}f\in L^p(\mR^d)\right\},
$$
where $\Delta^{\beta/2}$ is  defined by
$$
\Delta^{\beta/2}f:=\sF^{-1}\big(|\xi|^{\beta} \sF f\big),\quad \forall f\in \sS(\mR^d).
$$
We will need the following fact: for $\beta\in(0,1]$ and $p\in(d/\beta,\infty]$,
there is a constant $c=c(p,d,\beta)>0$ such that for all $f\in H^{\beta,p}(\mR^d)$,
\begin{align*}
\Big\|\sup_{y\not=0}|y|^{-\beta}|f(x+y)-f(x)|\Big\|_p\leq c\|f\|_{H^{\beta,p}}.
\end{align*}
Recall that the Hardy-Littlewood maximal function of a function $f$ is defined by
$$
\cM f(x):=\sup_{0<r<\infty}\frac{1}{|B_r|}\int_{B_r}|f(x+y)|\dif y,
$$
where $|B_r|$ denotes the Lebesgue measure of $B_r$. Then we have that for every $f\in H^{1,1}_{loc}(\mR^d)\cap C^1(\mR^d)$, there exists a constant $C_{d}>0$ such that
\begin{align}
|f(x)-f(y)|\leq C_{d}|x-y|\Big(\cM|\nabla f|(x)+\cM|\nabla f|(y)\Big),\quad \forall x,y\in\mR^d,   \label{w11}
\end{align}
and for $p\in(1,\infty]$, there exists a constant $C_{d,p}>0$ such that
\begin{align}
\|\cM f\|_p\leq C_{d,p}\|f\|_p.   \label{mf}
\end{align}
The following relationship between $B^\beta_{p,\infty}(\mR^d)$ and $H^{\beta,p}(\mR^d)$ can be found in \cite{Tri}: for $0<\beta\leq2$, $\eps\in(0,\beta)$ and $p\geq 2$, it holds that
\begin{align}
H^{\beta,p}(\mR^d)\subseteq B^{\beta}_{p,\infty}(\mR^d)\subseteq H^{\beta-\eps,p}(\mR^d).  \label{emb2}
\end{align}

Bernstein type inequalities are fundamental tools for studying differential equations using the Littlewood-Paley theory. We recall the following result, see \cite{B-C-D} or \cite{W-Z}.

\bl[Bernstein's inequality]
Let $1\leq p\leq q\leq \infty$. For any $0\leq k\in\mN$, there exists a constant $C_1=C_1(d,p,q,k,\alpha)>0$ such that for all $f\in\sS'(\mR^d)$ and $j\geq-1$,
\begin{align}
\|\nabla^k\Lambda_jf\|_q\leq C_12^{(k+d(\frac{1}{p}-\frac{1}{q}))j}\|\Lambda_jf\|_p,   \label{ber11}
\end{align}
and for $j\geq0$, $\alpha\in(0,2)$,
\begin{align*}
\|\Delta^{\alpha/2}\Lambda_jf\|_q\leq C_12^{(\alpha+d(\frac{1}{p}-\frac{1}{q}))j}\|\Lambda_jf\|_p.
\end{align*}
\el

The following commutator estimates can be found in \cite[Lemma 2.3]{C-Z-Z}.
\bl
Let $p,p_1,p_2,q_1,q_2\in[1,\infty]$ with $\frac{1}{p}=\frac{1}{p_1}+\frac{1}{p_2}$ and $\frac{1}{q_1}+\frac{1}{q_2}=1$. Then, for $\beta_1\in(0,1)$ and $\beta_2\in[-\beta_1,0]$, there exists a constant $C_2=C_2(d,p,p_1,p_2,\beta_1,\beta_2)>0$ such that for all $j\geq-1$,
\begin{equation}
\|[\Lambda_j,f]g\|_p\leq C_22^{-(\beta_1+\beta_2)j}\left\{\begin{aligned}
&\|f\|_{B^{\beta_1}_{p_1,\infty}}\|g\|_{p_2},\quad\quad\,\,\textrm{ if\,\,\, $\beta_2=0$};\\
&\|f\|_{B^{\beta_1}_{p_1,\infty}}\|g\|_{B^{\beta_2}_{p_2,\infty}},\quad\textrm{ if\,\,\, $\beta_2\neq0$\,\,\,and\,\,\, $\beta_1+\beta_2>0$};\\
&\|f\|_{B^{\beta_1}_{p_1,q_1}}\|g\|_{B^{\beta_2}_{p_2,q_2}},\quad\textrm{ if\,\,\,  $\beta_1+\beta_2=0$},\label{comm}
\end{aligned}
\right.
\end{equation}
where $[\Lambda_j,f]g:=\Lambda_j(fg)-f\Lambda_jg$.
\el

\section{Supercritical and critical parabolic equations}

In this section, we use the Littlewood-Paley theory to study the non-local partial differential equation corresponding to SDE (\ref{sdee}).
Below we fix  $T>0$. For $\lambda\geq0$, consider the following non-local parabolic PDE on $[0,T]\times\mR^d$:
\begin{equation} \label{pde}
\left\{ \begin{aligned}
&\p_tu(t,x)=\sL^{\sigma}_\nu u(t,x)-\lambda u(t,x)+b(t,x)\cdot\nabla u(t,x)+f(t,x),\\
&u(0,x)=0,
\end{aligned} \right.
\end{equation}
where $\sL^\sigma_\nu$ is defined in (\ref{ope}). Unless otherwise specified, we  assume {\bf (H$^{\nu}$)} holds true.
We will study the equation (\ref{pde}) both in Besov spaces and in H\"older spaces.
For simplicity, given a $\beta\in(0,1]$ and a function $f:[0,T]\times\mR^d\rightarrow\mR^d$, we will write
$$
\|f\|_{L^\infty_TC_b^\beta}:=\sup_{(t,x)\in[0,T]\times\mR^d}|f(t,x)|+\sup_{t\in[0,T]}[f(t)]_{\beta},
$$
where $[\cdot]_{\beta}$ denotes the H\"older semi-norm of a function defined by
$$
[f(t)]_{\beta}:=\sup_{x\neq y}\frac{|f(t,x)-f(t,y)|}{|x-y|^\beta}.
$$
For $p\in[1,\infty]$, we write $L^\infty_p(T):=L^\infty([0,T];L^p(\mR^d))$ with norm
$$
\|f\|_{L^\infty_p(T)}:=\sup_{t\in[0,T]}\|f(t,\cdot)\|_p.
$$
Given $\beta\in\mR$ and $p,q\in[1,\infty]$, we also use
$$
\|f\|_{L^\infty_TB^\beta_{p,q}}:=\sup_{t\in[0,T]}\|f(t,\cdot)\|_{B^\beta_{p,q}}
$$
to denote the norm of a function $f$ in $L^\infty\big([0,T];B^\beta_{p,q}(\mR^d)\big)$.

\subsection{Constant coefficient}

For simplicity, we first consider the case that $\sigma(t,x,z)=\sigma_0(t,z)$, that is,
the coefficient $\sigma$ is independent of the $x$-variable.
Throughout this subsection, we always assume that:

\begin{enumerate}
\item [{\bf (H$^{\sigma_0}$)}] There exist constants $\kappa_0,\kappa_1>0$ such that for all $(t,z)\in\mR_+\times\mR^d$,
\begin{align*}
\kappa_0\leq \sigma_0(t,z)\leq \kappa_1.
\end{align*}
\end{enumerate}
As in (\ref{ope}), we write
$$
\sL^{\sigma_0}_\nu f(t,x):=\int_{\mR^d}\Big[f(t,x+z)-f(t,x)
-1_{\{\alpha\geq 1\}}1_{\{|z|\leq 1\}}z\cdot\nabla_x f(t,x)\Big]\sigma_0(t,z)\nu(\dif z).
$$
We will study the following non-local parabolic PDE on $[0,T]\times\mR^d$:
\begin{equation} \label{pde1}
\left\{ \begin{aligned}
&\p_tu(t,x)=\sL^{\sigma_0}_\nu u(t,x)-\lambda u(t,x)+b(t,x)\cdot\nabla u(t,x)+f(t,x),\\
&u(0,x)=0,
\end{aligned} \right.
\end{equation}
where $\lambda\geq0$ is a constant.

We prepare some useful auxiliary results first. Note that by Fourier transform, we have  for each $t\geq0$,
$$
\sF(\sL^{\sigma_0}_\nu f)(\xi)=\psi^{\sigma_0}_\nu(t,\xi)\cdot\sF f(\xi),
$$
where the symbol $\psi^{\sigma_0}_\nu(t,\xi)$ is given by
$$
\psi_\nu^{\sigma_0}(t,\xi)=\int_{\mR^d}\big(\e^{i\xi\cdot z}-1
-1_{\{\alpha\geq 1\}}1_{\{|z|\leq 1\}}i\xi\cdot z\big)\sigma_0(t,z)\nu(\dif z).
$$
We have the following upper bound estimate for the symbol.

\bl \label{le4}
Under {\bf (H$^\nu$)} and {\bf (H$^{\sigma_0}$)}, there exist constants $C_0, C_1>0$ such that for every $t\geq0$,
\begin{align}
{\rm{Re}}(\psi_\nu^{\sigma_0}(t,\xi))\leq -C_0|\xi|^\alpha+C_1.   \label{sym}
\end{align}
\el
\begin{proof}
By assumption and (\ref{nu}), we have
\begin{align*}
{\rm{Re}}(\psi_\nu^{\sigma_0}(t,\xi))&=\int_{\mR^d}\big(\cos(\xi\cdot z)-1\big)\sigma_0(t,z)\nu(\dif z)\\
&\leq\kappa_0\int_{\mR^d}\big(\cos(\xi\cdot z)-1\big)\nu_1(\dif z)+
2\kappa_1\nu(B^c_1)\\
&=\kappa_0\left(\int_0^\infty\frac{(\cos r-1)\dif r}{r^{1+\alpha}}\right)\int_{\mS_{d-1}}\!|\xi\cdot\theta|^\alpha\Sigma(\dif \theta)
+2\kappa_1\nu(B^c_1)\\
&\leq-C_0|\xi|^\alpha+C_1.
\end{align*}
The proof is complete.
\end{proof}

With (\ref{sym}) in hand, we can prove the following Bernstein type inequality by using standard Fourier analysis method.

\bl\label{le42}
For any $p\geq2$, there exist positive constants $C_0, C_1$ such that for every $f\in\sS'(\mR^d)$ and $j\geq0$,
\begin{align}
\int_{\mR^d}|\Lambda_jf|^{p-2}\Lambda_jf(\sL^{\sigma_0}_\nu\Lambda_jf)\dif x\leq -C_02^{\alpha j}\|\Lambda_jf\|_p^p+C_1\|\Lambda_jf\|_p^p.   \label{ber2}
\end{align}
\el
\begin{proof}
We adopt the argument in \cite{C-Z-Z}. In fact, following the same procedure as in \cite[Lemma 3.1]{C-Z-Z}, we can show that for any $p\geq2$ and any smooth function $g$,
\begin{align*}
|g(x)|^{p/2-2}g(x)\sL^{\sigma_0}_\nu g(x)\leq\tfrac{2}{p}(\sL^{\sigma_0}_\nu|g|^{p/2})(x).
\end{align*}
As a result, we can get by Plancherel's theorem and the estimate (\ref{sym}) that
\begin{align*}
&\int_{\mR^d}|\Lambda_jf|^{p-2}\Lambda_jf(\sL^{\sigma_0}_\nu\Lambda_jf)\dif x\leq \frac{2}{p}\int_{\mR^d}|\Lambda_jf|^{p/2}\sL^{\sigma_0}_\nu|\Lambda_jf|^{p/2}\dif x\\
&\leq -\frac{2C_0}{p}\int_{\mR^d}\Big(|\xi|^{\alpha/2}\sF(|\Lambda_jf|^{p/2})(\xi)\Big)^2\dif \xi+\frac{2C_1}{p}\int_{\mR^d}|\Lambda_jf|^{p}\dif x\\
&= -\frac{2C_0}{p}\int_{\mR^d}\Big((-\Delta)^{\alpha/4}|\Lambda_jf|^{p/2}\Big)^2\dif x+\frac{2C_1}{p}\|\Lambda_jf\|_p^p.
\end{align*}
Recall that (see \cite[Proposition 3.1]{CMZ}) for every $j\geq0$,
$$
\|(-\Delta)^{\alpha/4}|\Lambda_jf|^{p/2}\|_2^2\geq c_02^{\alpha j}\|\Lambda_jf\|_p^p,
$$
where $c_0>0$ is independent of $j$. The desired estimate (\ref{ber2}) follows immediately.
\end{proof}

Now, we prove the following result for equation (\ref{pde1}) in Besov spaces.

\bt\label{drift2}

(i) Assume that $0<\alpha\leq1$, $b\in L^\infty\big([0,T];B^\beta_{p,\infty}(\mR^d)\big)$ with $\beta>1-\alpha$ and $d/(\alpha+\beta-1)\vee 2<p\leq\infty$. Then, for any $\lambda\geq0$ and  $f\in L^\infty\big([0,T];B^\gamma_{q,\infty}(\mR^d)\big)$ with $\gamma\in[0,\beta]$, $2\leq q\leq p$ and $q\neq\infty$, there exists a unique solution $u\in L^\infty\big([0,T];B^{\alpha+\gamma}_{q,\infty}(\mR^d)\big)$ to (\ref{pde1}) such that
\begin{align}
\|u\|_{L^\infty_TB^{\alpha+\gamma}_{q,\infty}}\leq C_1\|f\|_{L^\infty_TB^\gamma_{q,\infty}},   \label{ess1}
\end{align}
where $C_1=C(d,T,q,\alpha,\gamma,\|b\|_{L^\infty_TB^{\beta}_{p,\infty}})$ is a positive constant. Moreover, for any $\eta\in[0,\alpha+\gamma)$, we have
\begin{align}
\|u\|_{L^\infty_TB^{\eta}_{q,\infty}}\leq C_\lambda\|f\|_{L^\infty_TB^\gamma_{q,\infty}}  , \label{ess11}
\end{align}
where $C_\lambda$ is a positive constant satisfying $C_\lambda\rightarrow0$ as $\lambda\rightarrow\infty$.

\vspace{2mm}
\noindent (ii) Assume that $\alpha\in(0,1)$, $b\in L^\infty\big([0,T];C^{1-\alpha}_b(\mR^d)\big)$ with $\|b\|_{L^\infty_TC^{1-\alpha}_b}$ small. Then, for any $\lambda\geq0$ and $f\in L^\infty\big([0,T];B^\gamma_{q,\infty}(\mR^d)\big)$ with $\gamma\in(0,1-\alpha)$ and $2\leq q<\infty$, there exists a unique solution $u\in L^\infty\big([0,T];B^{\alpha+\gamma}_{q,\infty}(\mR^d)\big)$ to (\ref{pde1}) such that
\begin{align}
\|u\|_{L^\infty_TB^{\alpha+\gamma}_{q,\infty}}\leq C_2\|f\|_{L^\infty_TB^\gamma_{q,\infty}},    \label{ess2}
\end{align}
where $C_2=C(d,T,q,\alpha,\gamma,\|b\|_{L^\infty_TC^{1-\alpha}_b})$ is a positive
constant.  Moreover, for any $\eta\in[0,\alpha+\gamma)$, \eqref{ess11} holds with
a positive constant $C_\lambda$ satisfying $C_\lambda\rightarrow0$ as $\lambda\rightarrow\infty$.

\vspace{2mm}
\noindent (iii) Assume that $\alpha\in(0,1)$, $b\in L^\infty\big([0,T];B^{1-\alpha}_{\infty,1}(\mR^d)\big)$ with $\|b\|_{L^\infty_TB^{1-\alpha}_{\infty,1}}$ small. Then, for any $\lambda\geq0$ and $f\in L^\infty\big([0,T];B^0_{q,\infty}(\mR^d)\big)$ with $2\leq q<\infty$, there exists a unique solution $u\in L^\infty\big([0,T];B^{\alpha}_{q,\infty}(\mR^d)\big)$ to (\ref{pde1}) such that
\begin{align}
\|u\|_{L^\infty_TB^{\alpha}_{q,\infty}}\leq C_3\|f\|_{L^\infty_TB^0_{q,\infty}},    \label{ess6}
\end{align}
where $C_3=C(d,T,q,\alpha,\|b\|_{L^\infty_TB^{1-\alpha}_{\infty,1}})$ is a positive
constant. Moreover, for any $\eta\in[0,\alpha)$, \eqref{ess11} holds with $\gamma=0$ and
a positive constant $C_\lambda$ satisfying $C_\lambda\rightarrow0$ as $\lambda\rightarrow\infty$.
\et

\begin{proof}
It is well known that the non-local PDE (\ref{pde1}) has a unique smooth solution $u$ if
$$
b,f\in C_b^\infty(\mR_+\times\mR^{d}),
$$
see e.g., \cite{Zh}. Thus, we will focus on proving the a priori estimates (\ref{ess1}), (\ref{ess2}) and (\ref{ess6}), and  the estimate (\ref{ess11}). Then the conclusions follow by
a standard mollification method, see also \cite{C-Z-Z}.

\vspace{2mm}
\noindent (i) Using the operator $\Lambda_j$ to act on both sides of (\ref{pde1}) and by (\ref{con}), we have
\begin{align*}
\p_t\Lambda_ju&=\Lambda_j\sL^{\sigma_0}_\nu u+\Lambda_j(b\cdot\nabla u)+\Lambda_jf-\lambda\Lambda_ju\\
&=\sL^{\sigma_0}_\nu\Lambda_ju+[\Lambda_j,b\cdot\nabla]u+b\cdot\nabla\Lambda_ju+\Lambda_jf-\lambda\Lambda_ju.
\end{align*}
For $q\geq2$, multiplying both sides of the above equality by $|\Lambda_ju|^{q-2}\Lambda_ju$ and then integrating with respect to $x$ yield that
\begin{align*}
\tfrac{1}{q}\p_t\|\Lambda_ju\|_q^q&=\int_{\mR^d}(\sL^{\sigma_0}_\nu\Lambda_ju)|\Lambda_ju|^{q-2}\Lambda_ju\dif x+\int_{\mR^d}\big([\Lambda_j,b\cdot\nabla]u\big)|\Lambda_ju|^{q-2}\Lambda_ju\dif x\\
&\quad+\int_{\mR^d}(b\cdot\nabla\Lambda_ju)|\Lambda_ju|^{q-2}\Lambda_ju\dif x+\int_{\mR^d}(\Lambda_jf)|\Lambda_ju|^{q-2}\Lambda_ju\dif x-\lambda\|\Lambda_ju\|_q^q\\
&=:\cI_j^1+\cI^2_j+\cI^3_j+\cI^4_j-\lambda\|\Lambda_ju\|_q^q.
\end{align*}
For $\cI_j^1$, we have by Bernstein's inequalities (\ref{ber2}) that there exist
constants $\kappa, c_0>0$ such that for all $j\geq0$,
$$
\cI_j^1\leq -\kappa2^{\alpha j}\|\Lambda_ju\|_q^q+c_0\|\Lambda_ju\|_q^q.
$$
For $\cI_j^2$, we have by H\"older's inequality, (\ref{comm}) with $1/\hat p=1/q-1/p$ and the embedding (\ref{emb0}) and (\ref{pp}) that for every $\gamma\in[0,\beta]$, there exists a constant $c_1>0$ such that for all $j\geq-1$,
\begin{align*}
\cI_j^2\leq \|[\Lambda_j,b\cdot\nabla]u\|_q\|\Lambda_ju\|_q^{q-1}\leq c_12^{-\gamma j}\|b\|_{B^\beta_{p,\infty}}\|u\|_{B^{1-\beta+\gamma}_{\hat p,1}}\|\Lambda_ju\|_q^{q-1}.
\end{align*}
For the third term, we write
\begin{align*}
\cI^3_j&=\int_{\mR^d}\big((b-S_jb)\cdot\nabla\Lambda_ju)|\Lambda_ju|^{q-2}\Lambda_ju\dif x\\
&\quad+\int_{\mR^d}(S_jb\cdot\nabla\Lambda_ju)|\Lambda_ju|^{q-2}\Lambda_ju\dif x=:\cI^{31}_j+\cI^{32}_j.
\end{align*}
By H\"older's inequality and the definition of $S_j$, it is easy to see that
\begin{align*}
\cI^{31}_j\leq\sum_{k\geq j}\|\Lambda_kb\cdot\nabla\Lambda_ju\|_q\|\Lambda_ju\|_q^{q-1}.
\end{align*}
Using H\"older's inequality again with $1/\hat p=1/q-1/p$ and Bernstein's inequality (\ref{ber11}), we can deduce that for some constant $c_2>0$,
\begin{align*}
\|\Lambda_kb\cdot\nabla\Lambda_ju\|_q\leq \|\Lambda_kb\|_p\|\nabla\Lambda_ju\|_{\hat p}\leq c_2\|\Lambda_kb\|_p2^{(1+\frac{d}{p})j}\|\Lambda_ju\|_q.
\end{align*}
Thus we have
\begin{align*}
\cI^{31}_j&\leq\sum_{k\geq j}c_2\|\Lambda_kb\|_p\cdot2^{(1+\frac{d}{p})j}\|\Lambda_ju\|_q^{q}\leq\sum_{k\geq j}2^{-\beta k}\cdot c_22^{(1+\frac{d}{p})j}\|b\|_{B^\beta_{p,\infty}}\|\Lambda_ju\|_q^{q}\\
&\leq c_22^{(1+\frac{d}{p}-\beta)j}\|b\|_{B^\beta_{p,\infty}}\|\Lambda_ju\|_q^{q},\quad \forall j\geq-1.
\end{align*}
For $\cI^{32}_j$, we have by the divergence theorem and Bernstein's inequality (\ref{ber11}) that there exists a constant $c_3>0$ such that for all $j\geq-1$,
\begin{align*}
\cI^{32}_j&=\frac{1}{q}\int_{\mR^d}(S_jb\cdot\nabla)|\Lambda_ju|^{q}\dif x=-\frac{1}{q}\int_{\mR^d}(S_j\div b)|\Lambda_ju|^{q}\dif x\\
&\leq\frac{1}{q}\sum_{k\leq j-1}\|\Lambda_k\div b\|_\infty\|\Lambda_ju\|_q^q\leq \sum_{k\leq j-1}c_32^{(1+\frac{d}{p})k}\|\Lambda_kb\|_p\|\Lambda_ju\|_q^q\\
&\leq \sum_{k\leq j-1}c_32^{(1+\frac{d}{p}-\beta)k}\|b\|_{B^\beta_{p,\infty}}\|\Lambda_ju\|_q^q\leq c_32^{(1+\frac{d}{p}-\beta)j}\|b\|_{B^\beta_{p,\infty}}\|\Lambda_ju\|_q^{q}.
\end{align*}
For the last term, it is easy to see that
$$
\cI^4_j\leq \|\Lambda_jf\|_q\|\Lambda_ju\|_q^{q-1}.
$$
Combining the above estimates, we finally arrive at
\begin{align}
\tfrac{1}{q}\p_t\|\Lambda_ju\|_q^q&\leq -\kappa 1_{\{j\geq0\}}2^{\alpha j}\|\Lambda_ju\|_q^q-(\lambda-c_0)\|\Lambda_ju\|_q^q\no\\
&\quad+c_42^{(1+\frac{d}{p}-\beta)j}\|b\|_{B^\beta_{p,\infty}}\|\Lambda_ju\|_q^{q}\no\\
&\quad+c_42^{-\gamma j}\|b\|_{B^\beta_{p,\infty}}\|u\|_{B^{1-\beta+\gamma}_{\hat p,1}}\|\Lambda_ju\|_q^{q-1}+\|\Lambda_jf\|_q\|\Lambda_ju\|_q^{q-1}, \label{case}
\end{align}
where $c_4>0$ is a constant. By the assumption that
$$
1-\beta+d/p<\alpha\quad\text{and}\quad b\in L^\infty\big([0,T];B^\beta_{p,\infty}(\mR^d)\big),
$$
we can use Young's inequality and divide both sides of (\ref{case}) by $\|\Lambda_ju\|_q^{q-1}$ to get that for some $\kappa_0,\kappa_1>0$ and all $j\geq -1$,
\begin{align*}
\p_t\|\Lambda_ju\|_q\leq -(\kappa_02^{\alpha j}+\lambda-\kappa_1)\|\Lambda_ju\|_q
+c_52^{-\gamma j}\|u\|_{B^{1-\beta+\gamma}_{\hat p,1}}+\|\Lambda_jf\|_q,
\end{align*}
where $c_5>0$ depends on $\|b\|_{L^\infty_TB^\beta_{p,\infty}}$. This in turn implies by Gronwall's inequality that there exists a $c_6>0$ such that for all $j\geq -1$,
\begin{align}
\|\Lambda_ju(t)\|_q&\leq c_6\!\!\int_0^t\!\e^{-(\kappa_02^{\alpha j}+\lambda-\kappa_1)(t-s)}\Big(2^{-\gamma j}\|u(s)\|_{B^{1-\beta+\gamma}_{\hat p,1}}+\|\Lambda_jf(s)\|_q\Big)\dif s\no\\
&\leq c_62^{-\gamma j}\!\!\int_0^t\!\e^{-(\kappa_02^{\alpha j}+\lambda-\kappa_1)(t-s)}\Big(\|u(s)\|_{B^{1-\beta+\gamma}_{\hat p,1}}+\|f(s)\|_{B^\gamma_{q,\infty}}\Big)\dif s.  \label{ggg}
\end{align}
Multiplying both sides by $2^{(\alpha+\gamma)j}$ and taking supremum over $j$, we get that for a constant $c_7>0$, it holds that for every $t\in(0,T]$,
\begin{align}
\|u\|_{L^\infty_tB^{\alpha+\gamma}_{q,\infty}}&\leq c_7\Big(\|u\|_{L^\infty_tB^{1-\beta+\gamma}_{\hat p,1}}+\|f\|_{L^\infty_tB^\gamma_{q,\infty}}\Big).  \label{aa}
\end{align}
Notice that since $p>d/(\alpha+\beta-1)$, we have by (\ref{emb}) that for $\theta\in(0,\alpha+\beta-1-d/p)$,
$$
B^{\alpha-\theta+\gamma}_{q,\infty}(\mR^d)\subseteq B^{1-\beta+\gamma}_{\hat p,1}(\mR^d).
$$
Thus by  (\ref{int}) and (\ref{ineq}), we have that, for every $\eps>0$, there exists a $c_\eps>0$ such that
$$
\|u(t)\|_{B^{1-\beta+\gamma}_{\hat p,1}}\leq\eps\|u(t)\|_{B^{\alpha+\gamma}_{q,\infty}}+c_\eps\|u(t)\|_{B^{\gamma}_{q,\infty}}.
$$
Plugging this back into (\ref{aa}) and choosing $\eps$ small enough, we get that for every $t\in(0,T]$,
\begin{align}
\|u\|_{L^\infty_tB^{\alpha+\gamma}_{q,\infty}}\leq c_8\big(\|u\|_{L^\infty_tB^{\gamma}_{q,\infty}}+\|f\|_{L^\infty_TB^\gamma_{q,\infty}}\big), \label{77}
\end{align}
where $c_8>0$ is a constant. On the other hand, using (\ref{77}) and  (\ref{ggg}), we  also have
\begin{align*}
\|u\|_{L^\infty_tB^\gamma_{q,\infty}}&\leq c_9\!\!\int_0^t\!\e^{-(\lambda-\kappa_1)(t-s)}\|u\|_{L^\infty_sB^{1-\beta+\gamma}_{\hat p,1}}\dif s+c_\lambda\|f\|_{L^\infty_TB^\gamma_{q,\infty}}\\
&\leq c_9\!\!\int_0^t\!\|u\|_{L^\infty_sB^{\gamma}_{q,\infty}}\dif s+c_\lambda\|f\|_{L^\infty_TB^\gamma_{q,\infty}},
\end{align*}
where $c_9$ and $c_\lambda$ are  positive constants with $c_\lambda$
satisfying $c_\lambda\rightarrow0$ as $\lambda\rightarrow\infty$. Gronwall's inequality yields that
\begin{align*}
\|u\|_{L^\infty_TB^\gamma_{q,\infty}}\leq c_\lambda\|f\|_{L^\infty_TB^\gamma_{q,\infty}}.
\end{align*}
This together with (\ref{77}) implies (\ref{ess1}), and (\ref{ess11}) follows by interpolation.

\vspace{2mm}
(ii) Recall that by (\ref{bh}), we have $B^{1-\alpha}_{\infty,\infty}(\mR^d)=C_b^{1-\alpha}(\mR^d)$. Since in case (i), we can take $p=\infty$. Thus one can check that all the estimates for $\cI_j^1$, $\cI_j^3$ and $\cI_j^4$ in the proof of i) hold with $p=\infty$ and $\beta=1-\alpha$, i.e., there exist constants $\kappa, c_0, c_1>0$ such that
\begin{align*}
\cI_j^1+\cI_j^3+\cI_j^4\leq&-\kappa 1_{\{j\geq0\}}2^{\alpha j}\|\Lambda_ju\|_q^q+c_0\|\Lambda_ju\|_q^q\\
&+c_12^{\alpha j}\|b\|_{C^{1-\alpha}_b}\|\Lambda_ju\|_q^{q}+\|\Lambda_jf\|_q\|\Lambda_ju\|_q^{q-1}.
\end{align*}
For $\cI_j^2$, since  $\gamma\in(0,1-\alpha)$, we can use (\ref{comm}) to get that for some constant $c_2>0$,
\begin{align*}
\cI_j^2\leq c_22^{-\gamma j}\|b\|_{C^{1-\alpha}_b}\|u\|_{B^{\alpha+\gamma}_{q,\infty}}\|\Lambda_ju\|_q^{q-1}.
\end{align*}
Using the same procedures as in the proof of (i), we can arrive at
\begin{align}
\tfrac{1}{q}\p_t\|\Lambda_ju\|_q&\leq -\kappa1_{\{j\geq0\}}2^{\alpha j}\|\Lambda_ju\|_q-(\lambda-c_0)\|\Lambda_ju\|_q\no\\
&\quad+c_12^{\alpha j}\|b\|_{C^{1-\alpha}_b}\|\Lambda_ju\|_q+c_22^{-\gamma j}\|b\|_{C^{1-\alpha}_b}\|u\|_{B^{\alpha+\gamma}_{q,\infty}}+\|\Lambda_jf\|_q. \label{small}
\end{align}
Thus, if $\|b\|_{L^\infty_TC^{1-\alpha}_b}$ is small enough so that
$$
c_1\|b\|_{L^\infty_TC^{1-\alpha}_b}<\kappa,
$$
then we can get by Gronwall's inequality that for all $j\geq -1$, there exist constants $\kappa_0, \kappa_1, c_3, c_4>0$ such that
\begin{align}
\|\Lambda_ju(t)\|_q&\leq c_32^{-\gamma j}\!\!\int_0^t\!\e^{-(\kappa_02^{\alpha j}+\lambda-\kappa_1)(t-s)}\dif s\Big(\|b\|_{L^\infty_TC^{1-\alpha}_b}\|u\|_{L^\infty_TB^{\alpha+\gamma}_{q,\infty}}+\|f\|_{L^\infty_TB^\gamma_{q,\infty}}\Big) \label{333}\\
&\leq c_42^{-(\alpha+\gamma)j}\Big(\|b\|_{L^\infty_TC^{1-\alpha}_b}\|u\|_{L^\infty_TB^{\alpha+\gamma}_{q,\infty}}+\|f\|_{L^\infty_TB^\gamma_{q,\infty}}\Big).\no
\end{align}
This in particular implies that
\begin{align*}
\|u\|_{L^\infty_TB^{\alpha+\gamma}_{q,\infty}}\leq c_4\Big(\|b\|_{L^\infty_TC^{1-\alpha}_b}\|u\|_{L^\infty_TB^{\alpha+\gamma}_{q,\infty}}+\|f\|_{L^\infty_TB^\gamma_{q,\infty}}\Big).
\end{align*}
Now we further take $\|b\|_{L^\infty_TC^{1-\alpha}_b}$ small enough so that
$$
c_4\\|b\|_{L^\infty_TC^{1-\alpha}_b}<1.
$$
This in turn yields that
\begin{align*}
\|u\|_{L^\infty_TB^{\alpha+\gamma}_{q,\infty}}\leq c_5\|f\|_{L^\infty_TB^\gamma_{q,\infty}},
\end{align*}
where $c_5>0$ is a constant. Thus (\ref{ess2}) is true. Plugging this back into (\ref{333}) we get that
\begin{align*}
\|u\|_{L^\infty_TB^{\gamma}_{q,\infty}}\leq c_\lambda\|f\|_{L^\infty_TB^\gamma_{q,\infty}},
\end{align*}
where $c_\lambda\rightarrow0$ as $\lambda\rightarrow\infty$. This together with interpolation implies (\ref{ess11}).

\vspace{2mm}
(iii) We only give the main difference with the proof of (ii). In this case, since $\gamma=0$, we can take $q_1=1$ and $q_2=\infty$ in (\ref{comm}) to get that for some $c_1>0$,
\begin{align*}
\cI_j^2\leq c_1\|b\|_{B^\beta_{\infty,1}}\|u\|_{B^{\alpha}_{q,\infty}}\|\Lambda_ju\|_q^{q-1}.
\end{align*}
Thus we can get
\begin{align*}
\tfrac{1}{q}\p_t\|\Lambda_ju\|_q&\leq -\kappa1_{\{j\geq0\}}2^{\alpha j}\|\Lambda_ju\|_q-(\lambda-c_0)\|\Lambda_ju\|_q\\
&\quad+c_12^{\alpha j}\|b\|_{C^{1-\alpha}_b}\|\Lambda_ju\|_q+c_2\|b\|_{B^{1-\alpha}_{\infty,1}}\|u\|_{B^{\alpha}_{q,\infty}}+\|\Lambda_jf\|_q,
\end{align*}
where $c_0, c_2$ are as in part (ii) and $\kappa$ is as in part (i).
Notice that
$$
\|b\|_{C^{1-\alpha}_b}\leq \|b\|_{B^{1-\alpha}_{\infty,1}}.
$$
Following the same arguments as in the proof of ii), we can get the desired result.
The proof is complete.
\end{proof}

\br \label{ppp}
Inspecting the proofs of Lemmas \ref{le4}, \ref{le42}, Theorem \ref{drift2}.(ii) and \ref{drift2}.(iii) (particularly \eqref{small}), we can see that if the constant $\kappa_0$ in  {\bf (H$^\sigma_0$)} is large enough, then the smallness condition on $b$ in Theorem \ref{drift2}.(ii) and \ref{drift2}.(iii) can be dropped.
\er

\subsection{Variable coefficient}
In this subsection, we consider the non-local equation (\ref{pde}) in the variable coefficient case. Throughout this subsection, we assume that {\bf (H$^\sigma_1$)} holds true.
This in particular implies that
$$
\sigma\in L^\infty([0,T]\times\mR^d_z;C_b^\vartheta(\mR^d_x)),
$$
where $\vartheta$ is the constant in {\bf (H$^\sigma_1$)} .
For simplify,  we will denote by $\|\sigma\|_{L^\infty_\infty(T)C_b^\vartheta}$ the norm of $\sigma$ in $L^\infty([0,T]\times\mR^d_z;C_b^\vartheta(\mR^d_x))$.

Fix $z\in\mR^d$ below. Given a function $f$ on $\mR^d$, we introduce the shift operator
\begin{align}
\sT_zf(x):=f(x+z)-f(x).    \label{sh}
\end{align}
Define the commutator
$$
[\Lambda_j,\sL^\sigma_\nu]f:=\Lambda_j(\sL^\sigma_\nu f)-\sL^\sigma_\nu(\Lambda_jf).
$$
We first establish the following commutator estimate.

\bl\label{newc}
For any $1<p\leq\infty$, $\bar\vartheta>0$ and $\bar\vartheta-\vartheta<\gamma\leq\bar\vartheta$, there exists a constant $C_1=C_1(d,p,\vartheta,\gamma)>0$ such that for any $f\in\sS'(\mR^d)$ and $j\geq-1$,
\begin{align*}
\|[\Lambda_j,\sL^\sigma_\nu]f\|_p\leq C_12^{-(\vartheta-\bar\vartheta+\gamma) j}\|\sigma\|_{L^\infty_\infty(T)C^{\vartheta}_b}\|f\|_{B^{\alpha-\bar\vartheta+\gamma}_{p,\infty}}.
\end{align*}
\el
\begin{proof}
By definition, we can write
$$
\sL^{\sigma}_\nu f(x)=\int_{\mR^d}\sT_zf(x)\sigma(t,x,z)\nu(\dif z).
$$
Using the Bony decomposition \eqref{bon}, we have
\begin{align*}
\sT_zf(x)\sigma(t,x,z)=\cT_{\sigma(t,x,z)}\sT_zf(x)+\cT_{\sT_zf(x)}\sigma(t,x,z)+\cR\big(\sigma(t,x,z),\sT_zf(x)\big).
\end{align*}
Thus,
\begin{align*}
\Lambda_j(\sL^\sigma_\nu f)(x)&=\int_{\mR^d}\Big[\Lambda_j\Big(\cT_{\sigma(t,x,z)}\sT_zf(x)\Big)+\Lambda_j\Big(\cT_{\sT_zf(x)}\sigma(t,x,z)\Big)\\
&\quad+\Lambda_j\Big(\cR\big(\sigma(t,x,z),\sT_zf(x)\big)\Big)\Big]\nu(\dif z).
\end{align*}
Similarly, we can write
\begin{align*}
\sL^\sigma_\nu(\Lambda_jf)(x)&=\int_{\mR^d}\Lambda_j\sT_zf(x)\cdot\sigma(t,x,z)\nu(\dif z)\\
&=\int_{\mR^d}\Big[\cT_{\sigma(t,x,z)}\Lambda_j\sT_zf(x)+\cT_{\Lambda_j\sT_zf(x)}\sigma(t,x,z)\\
&\qquad\qquad\qquad\qquad+\cR\big(\sigma(t,x,z),\Lambda_j\sT_zf(x)\big)\Big]\nu(\dif z).
\end{align*}
As a result, we have
\begin{align*}
[\Lambda_j,\sL^\sigma_\nu]f(x)&=\int_{\mR^d}[\Lambda_j,\cT_{\sigma(t,x,z)}]\sT_zf(x)\nu(\dif z)+\int_{\mR^d}\Lambda_j\Big(\cT_{\sT_zf(x)}\sigma(t,x,z)\Big)\nu(\dif z)\\
&\quad-\int_{\mR^d}\cT_{\Lambda_j\sT_zf(x)}\sigma(t,x,z)\nu(\dif z)+\int_{\mR^d}\Lambda_j\Big(\cR\big(\sigma(t,x,z),\sT_zf(x)\big)\Big)\nu(\dif z)\\
&\quad-\int_{\mR^d}\cR\big(\sigma(t,x,z),\Lambda_j\sT_zf(x)\big)\nu(\dif z)=:\cQ^1_j+\cQ^2_j+\cQ^3_j+\cQ^4_j+\cQ^5_j.
\end{align*}
Below, we will omit the arguments of the functions, and proceed to control each term. For $\cQ_j^1$, thanks to (\ref{00}), we can write
\begin{align*}
[\Lambda_j,\cT_{\sigma}]\sT_zf&=\sum_{|k-j|\leq 4}\big(\Lambda_j(S_{k-1}\sigma\cdot\Lambda_k\sT_zf)-S_{k-1}\sigma\cdot\Lambda_j\Lambda_k\sT_zf\big)\\
&=\sum_{|k-j|\leq 4}[\Lambda_j,S_{k-1}\sigma]\Lambda_k\sT_zf.
\end{align*}
Note that
\begin{align*}
\big|[\Lambda_j,S_{k-1}\sigma]\Lambda_k\sT_zf\big|&\leq\int_{\mR^d}|h_j(y)|\big|S_{k-1}\sigma(t,x-y,z)-S_{k-1}\sigma(t,x,z)\big|\cdot|\Lambda_k\sT_zf(x-y)|\dif y\\
&\leq\|\sigma\|_{L^\infty_\infty(T)C^{\vartheta}_b}\int_{\mR^d}|h_j(y)||y|^\vartheta\cdot|\Lambda_k\sT_zf(x-y)|\dif y.
\end{align*}
Hence, we have
\begin{align*}
\|\cQ_j^1\|_p&\leq\sum_{|k-j|\leq 4}\left\|\int_{\mR^d}\big|[\Lambda_j,S_{k-1}\sigma]\Lambda_k\sT_zf\big|\nu(\dif z)\right\|_p\\
&\leq \|\sigma\|_{L^\infty_\infty(T)C^{\vartheta}_b}\sum_{|k-j|\leq 4}\int_{\mR^d}|h_j(y)||y|^\vartheta\dif y\cdot\left\|\int_{\mR^d}|\Lambda_k\sT_zf(x)|\nu(\dif z)\right\|_p\\
&\preceq2^{-\vartheta j}\|\sigma\|_{L^\infty_\infty(T)C^{\vartheta}_b}\sum_{|k-j|\leq4}\left\|\int_{\mR^d}|\Lambda_k\sT_zf(x)|\nu(\dif z)\right\|_p.
\end{align*}
We write
\begin{align*}
\left\|\int_{\mR^d}|\Lambda_k\sT_zf(x)|\nu(\dif z)\right\|_p=\left\|\left(\int_{|z|\leq2^{-k}}+\int_{|z|>2^{-k}}\right)|\Lambda_k\sT_zf(x)|\nu(\dif z)\right\|_p=:I_1+I_2.
\end{align*}
For $I_1$, we have by the mean value theorem and  Bernstein's inequality (\ref{ber11}) that for a constant $\theta\in[0,1]$,
\begin{align*}
I_1&=\left\|\int_{|z|\leq2^{-k}}|z|\cdot|\nabla\Lambda_kf(x+\theta z)|\nu(\dif z)\right\|_p\\
&\leq\int_{|z|\leq2^{-k}}|z|\nu(\dif z)\cdot\|\nabla\Lambda_kf\|_p\preceq 2^{\alpha k}\|\Lambda_kf\|_p.
\end{align*}
For $I_2$, it is easy to see that
\begin{align*}
I_2\leq\int_{|z|>2^{-k}}\nu(\dif z)\|\Lambda_kf\|_p\preceq 2^{\alpha k}\|\Lambda_kf\|_p.
\end{align*}
As a result, we have
\begin{align*}
\|\cQ_j^1\|_p&\preceq2^{-\vartheta j}\|\sigma\|_{L^\infty_\infty(T)C^{\vartheta}_b}\sum_{|k-j|\leq4}2^{\alpha k}\|\Lambda_kf\|_{p}\\
&\leq2^{-\vartheta j}\|\sigma\|_{L^\infty_\infty(T)C^{\vartheta}_b}\|f\|_{B^{\alpha-\bar\vartheta+\gamma}_{p,\infty}}\sum_{|k-j|\leq4}2^{(\bar\vartheta-\gamma) k}\\
&\preceq2^{-(\vartheta-\bar\vartheta+\gamma) j}\|\sigma\|_{L^\infty_\infty(T)C^{\vartheta}_b}\|f\|_{B^{\alpha-\bar\vartheta+\gamma}_{p,\infty}}.
\end{align*}
Similarly, we can write
$$
\Lambda_j(\cT_{\sT_zf}\sigma)=\sum_{|k-j|\leq 4}\Lambda_j\big(S_{k-1}\sT_zf\cdot\Lambda_k\sigma\big).
$$
Hence,  we can control the second term by
\begin{align*}
\|\cQ_j^2\|_p&\leq\sum_{|k-j|\leq4}\sum_{m\leq k-2}\left\|\int_{\mR^d}|\Lambda_m\sT_zf\cdot\Lambda_k\sigma|\nu(\dif z)\right\|_p\\
&\leq\sum_{|k-j|\leq4}\sum_{m\leq k-2}\|\Lambda_k\sigma\|_\infty\left\|\int_{\mR^d}|\Lambda_m\sT_zf|\nu(\dif z)\right\|_p\\
&\preceq\|\sigma\|_{L^\infty_\infty(T)C^{\vartheta}_b}\|f\|_{B^{\alpha-\bar\vartheta+\gamma}_{p,\infty}}\sum_{|k-j|\leq4}2^{-\vartheta k}\sum_{m\leq k-2}2^{(\bar\vartheta-\gamma)m}\\
&\preceq2^{-(\vartheta-\bar\vartheta+\gamma) j}\|\sigma\|_{L^\infty_\infty(T)C^{\vartheta}_b}\|f\|_{B^{\alpha-\bar\vartheta+\gamma}_{p,\infty}}.
\end{align*}
For $\cQ_j^3$, we have
\begin{align*}
\|\cQ_j^3\|_p&\leq\sum_{k\geq j-2}\left\|\int_{\mR^d}|S_{k-1}\Lambda_j\sT_zf\cdot\Lambda_k\sigma|\nu(\dif z)\right\|_p\\
&\leq\sum_{k\geq j-2}\|\Lambda_k\sigma\|_\infty\left\|\int_{\mR^d}|\Lambda_j\sT_zf|\nu(\dif z)\right\|_p\\
&\preceq\|\sigma\|_{L^\infty_\infty(T)C^{\vartheta}_b}\sum_{k\geq j-2}2^{-\vartheta k}\left\|\int_{\mR^d}|\Lambda_j\sT_zf|\nu(\dif z)\right\|_p\\
&\preceq2^{-(\vartheta-\bar\vartheta+\alpha+\gamma) j}\|\sigma\|_{L^\infty_\infty(T)C^{\vartheta}_b}\|f\|_{B^{\alpha-\bar\vartheta+\gamma}_{p,\infty}}.
\end{align*}
Finally, since $\vartheta-\bar\vartheta+\gamma>0$, we have
\begin{align*}
\|\cQ_j^4\|_p&\leq\left\|\int_{\mR^d}|\Lambda_j(\cR(\sigma,\sT_zf))|\nu(\dif z)\right\|_p\leq\sum_{|i|\leq 1,k\geq j-4}\left\|\int_{\mR^d}|\Lambda_j(\Lambda_k\sigma\cdot\Lambda_{k-i}\sT_zf)|\nu(\dif z)\right\|_p\\
&\preceq\|\sigma\|_{L^\infty_\infty(T)C^{\vartheta}_b}\sum_{|i|\leq 1,k\geq j-4}2^{-\vartheta k}\left\|\int_{\mR^d}|\Lambda_{k-i}\sT_zf|\nu(\dif z)\right\|_p\\
&\preceq\|\sigma\|_{L^\infty_\infty(T)C^{\vartheta}_b}\|f\|_{B^{\alpha-\bar\vartheta+\gamma}_{p,\infty}}\sum_{|i|\leq 1,k\geq j-4}2^{-\vartheta k}2^{(\bar\vartheta-\gamma)(k-i)}\\
&\preceq2^{-(\vartheta-\bar\vartheta+\gamma) j}\|\sigma\|_{L^\infty_\infty(T)C^{\vartheta}_b}\|f\|_{B^{\alpha-\bar\vartheta+\gamma}_{p,\infty}},
\end{align*}
and similarly for $\cQ_j^5$, we have
\begin{align*}
\|\cQ_j^5\|_p&\leq\left\|\int_{\mR^d}|\cR(\sigma,\Lambda_j\sT_zf)|\nu(\dif z)\right\|_p\leq\sum_{|i|\leq1,|k-j|\leq1}\left\|\int_{\mR^d}|\Lambda_{k-i}\sigma\cdot\Lambda_k\Lambda_j\sT_zf|\nu(\dif z)\right\|_p\\
&\leq2^{-(\vartheta-\bar\vartheta+\gamma) j}\|\sigma\|_{L^\infty_\infty(T)C^{\vartheta}_b}\|f\|_{B^{\alpha-\bar\vartheta+\gamma}_{p,\infty}}.
\end{align*}
Combining the above estimates, we get the desired result.
\end{proof}

To study equation (\ref{pde}) with variable coefficients, we need to use the freezing
coefficient method.
To this end, we introduce the following freezing function: let $q\geq1$ and let $\phi\in C_0^\infty(\mR^d)$ be a non-negative function with support in the unit ball and satisfying
$\int_{\mR^d}\phi^q(x)\dif x=1$.
For $y\in\mR^d$ and $\delta\in(0,1)$, define
\begin{align}
\phi^\delta_y(x):=\delta^{-d/q}\phi(\delta^{-1}(x-y)).    \label{delta}
\end{align}
We prove the following result for equation (\ref{pde}) in Besov spaces.

\bt\label{drift3}
(i) Assume that $0<\alpha\leq1$, {\bf (H$^\sigma_1$)} holds, $b\in L^\infty\big([0,T];B^\beta_{p,\infty}(\mR^d)\big)$ with $\beta>1-\alpha$ and $d/(\alpha+\beta-1)\vee 2<p\leq\infty$. Then, for any $\lambda\geq0$ and  $f\in L^\infty\big([0,T];B^\gamma_{q,\infty}(\mR^d)\big)$ with $\gamma\in[0,\beta\wedge\vartheta]$, $2\leq q\leq p$ and $q\neq\infty$, there exists a unique solution $u\in L^\infty\big([0,T];B^{\alpha+\gamma}_{q,\infty}(\mR^d)\big)$ to (\ref{pde}) such that
\begin{align}
\|u\|_{L^\infty_TB^{\alpha+\gamma}_{q,\infty}}\leq C_1\|f\|_{L^\infty_TB^\gamma_{q,\infty}},   \label{es1}
\end{align}
where $C_1=C(d,T,q,\alpha,\gamma,\|\sigma\|_{L^\infty_\infty C_b^\vartheta},\|b\|_{L^\infty_TB^{\beta}_{p,\infty}})$ is a positive constant. Moreover, for any $\eta\in[0,\alpha+\gamma)$,
\begin{align}
\|u\|_{L^\infty_TB^{\eta}_{q,\infty}}\leq C_\lambda\|f\|_{L^\infty_TB^\gamma_{q,\infty}},   \label{ess111}
\end{align}
where $C_\lambda$ is a positive constant satisfying $C_\lambda\rightarrow0$ as $\lambda\rightarrow\infty$.

\vspace{2mm}
\noindent (ii) Assume that $\alpha\in(0,1)$, {\bf (H$^\sigma_1$)} holds, $b\in L^\infty\big([0,T];C^{1-\alpha}_b(\mR^d)\big)$ with $\|b\|_{L^\infty_TC^{1-\alpha}_b}$ small. Then, for any $\lambda\geq0$ and $f\in L^\infty\big([0,T];B^\gamma_{q,\infty}(\mR^d)\big)$ with $\gamma\in(0,(1-\alpha)\wedge\vartheta)$ and $2\leq q<\infty$, there exists a unique solution $u\in L^\infty\big([0,T];B^{\alpha+\gamma}_{q,\infty}(\mR^d)\big)$ to (\ref{pde}) with
\begin{align}
\|u\|_{L^\infty_TB^{\alpha+\gamma}_{q,\infty}}\leq C_2\|f\|_{L^\infty_TB^\gamma_{q,\infty}},    \label{ess22}
\end{align}
where $C_2=C(d,T,q,\alpha,\gamma,\|\sigma\|_{L^\infty_\infty C_b^\vartheta},\|b\|_{L^\infty_TC^{1-\alpha}_b})$ is a positive
constant. Moreover, for any $\eta\in[0,\alpha+\gamma)$, (\ref{ess111}) holds with
a positive constant $C_\lambda$ satisfying $C_\lambda\rightarrow0$ as $\lambda\rightarrow\infty$.

\vspace{2mm}
\noindent (iii) Assume that $\alpha\in(0,1)$, {\bf (H$^\sigma_1$)} holds, $b\in L^\infty\big([0,T];B^{1-\alpha}_{\infty,1}(\mR^d)\big)$ with $\|b\|_{L^\infty_TB^{1-\alpha}_{\infty,1}}$ small. Then, for any $\lambda\geq0$ and $f\in L^\infty\big([0,T];B^0_{q,\infty}(\mR^d)\big)$ with $2\leq q<\infty$, there exists a unique solution $u\in L^\infty\big([0,T];B^{\alpha}_{q,\infty}(\mR^d)\big)$ to (\ref{pde}) with
\begin{align}
\|u\|_{L^\infty_TB^{\alpha}_{q,\infty}}\leq C_3\|f\|_{L^\infty_TB^0_{q,\infty}},    \label{ess66}
\end{align}
where $C_3=C(d,T,q,\alpha,\|\sigma\|_{L^\infty_\infty C_b^\vartheta},\|b\|_{L^\infty_TB^{1-\alpha}_{\infty,1}})$ is a positive
constant. Moreover, for any $\eta\in[0,\alpha)$, (\ref{ess111}) holds with $\gamma=0$ and
a positive constant $C_\lambda$ satisfying $C_\lambda\rightarrow0$ as $\lambda\rightarrow\infty$.
\et
\begin{proof}
By the classical continuity method, it suffices to prove the a priori estimate (\ref{es1}), (\ref{ess22}) and (\ref{ess66}), and the estimate (\ref{ess111}).

\vspace{2mm}
(i)
Using the operator $\Lambda_j$ to act on both sides of (\ref{pde}), we get
\begin{align*}
\p_t\Lambda_ju&=\Lambda_j(\sL^{\sigma}_\nu u)+\Lambda_j(b\cdot\nabla u)-\lambda\Lambda_ju+\Lambda_jf\\
&=\sL^{\sigma_0}_\nu\Lambda_ju+b\cdot\nabla\Lambda_ju-\lambda\Lambda_ju\\
&\quad+\Lambda_jf+[\Lambda_j,\sL^\sigma_\nu]u+(\sL^\sigma_\nu-\sL^{\sigma_0}_\nu)\Lambda_ju+[\Lambda_j,b\cdot\nabla]u,
\end{align*}
where $\sigma_0(t, z):=\sigma(t,y,z)$ with $y\in\mR^d$ being fixed.
Let $\phi_y^\delta$ be given by (\ref{delta}) with $q$ being in the statement of the theorem,
multiplying the above equation by $\phi_y^\delta$, we obtain the following
equation with constant coefficients:
\begin{align*}
\p_t(\Lambda_ju\cdot\phi_y^\delta)&=\sL^{\sigma_0}_\nu(\Lambda_ju\cdot\phi_y^\delta)+b\cdot\nabla(\Lambda_ju\cdot\phi_y^\delta) -\lambda(\Lambda_ju\cdot\phi_y^\delta)\\
&\quad+\Lambda_jf\cdot\phi_y^\delta+\big(\sL^{\sigma_0}_\nu\Lambda_ju\cdot\phi_y^\delta-\sL^{\sigma_0}_\nu(\Lambda_ju\cdot\phi_y^\delta)\big)-b\cdot \Lambda_ju\cdot\nabla\phi_y^\delta\\
&\quad+[\Lambda_j,\sL^\sigma_\nu]u\cdot\phi_y^\delta+(\sL^\sigma_\nu-\sL^{\sigma_0}_\nu)\Lambda_ju\cdot\phi_y^\delta +[\Lambda_j,b\cdot\nabla]u\cdot\phi_y^\delta.
\end{align*}
For simplicity, we define
\begin{align*}
\tilde f^\delta_j(t,x,y)&:=\Lambda_jf\cdot\phi_y^\delta-b\cdot \Lambda_ju\cdot\nabla\phi_y^\delta+[\Lambda_j,b\cdot\nabla]u\cdot\phi_y^\delta+[\Lambda_j,\sL^\sigma_\nu]u\cdot\phi_y^\delta\\
&\quad+(\sL^\sigma_\nu-\sL^{\sigma_0}_\nu)\Lambda_ju\cdot\phi_y^\delta +\big(\sL^{\sigma_0}_\nu\Lambda_ju\cdot\phi_y^\delta -\sL^{\sigma_0}_\nu(\Lambda_ju\cdot\phi_y^\delta)\big).
\end{align*}
Then, repeating the argument in the proof of Theorem \ref{drift2}.(i),
we get that for any $2\leq q<\infty$ and $j\geq-1$, there exist constants $\kappa_0, \kappa_1>0$ such that
\begin{align*}
\|\Lambda_ju(t)\cdot\phi_y^\delta\|_q\preceq \int_0^t-(\kappa_02^{\alpha j}+\lambda-\kappa_1)\|\Lambda_ju(s)\cdot\phi_y^\delta\|_q\dif s
+\int_0^t\|\tilde f_j^\delta(s,\cdot,y)\|_q\dif s.
\end{align*}
Taking $L^q$-norm with respect to the $y$ variable on the left hand side of the above inequality and noticing that $\|\phi_y^\delta\|_q=1$,
we have
$$
\left(\int_{\mR^d}\|\Lambda_ju(t)\cdot\phi_y^\delta\|_q^q\dif y\right)^{1/q}=\left(\int_{\mR^d}|\Lambda_ju(t)|^q\int_{\mR^d}|\phi_y^\delta|^q\dif y\dif x\right)^{1/q}=\|\Lambda_ju\|_q.
$$
Thus we have
\begin{align*}
\|\Lambda_ju(t)\|_q\preceq\int_0^t-(\kappa_02^{\alpha j}+\lambda-\kappa_1)\|\Lambda_ju(s)\|_q\dif s+\int_0^t\left(\int_{\mR^d}\|\tilde f_j^\delta(s,\cdot,y)\|_q^q\dif y\right)^{1/q}\dif s.
\end{align*}
Below, we will omit the $t$-variable and proceed to show that there exist constants $\kappa_2,\eps,C_\eps>0$ such that
\begin{align}
\left(\int_{\mR^d}\|\tilde f_j^\delta(\cdot,y)\|_q^q\dif y\right)^{1/q}\leq \kappa_22^{-\gamma j}\Big((\eps+\delta^\vartheta)\|u\|_{B^{\alpha+\gamma}_{q,\infty}}+C_\eps\|u\|_{B^{\gamma}_{q,\infty}}+\|f\|_{B^{\gamma}_{q,\infty}}\Big). \label{ess}
\end{align}
In fact, we can write
$$
\left(\int_{\mR^d}\|\tilde f_j^\delta(\cdot,y)\|_q^q\dif y\right)^{1/q}\leq \cR^1_j+\cR_j^2+\cR_j^3+\cR_j^4+\cR_j^5+\cR_j^6.
$$
For the first term, it is easy to see that
\begin{align*}
\cR_j^1=\left(\int_{\mR^d}\|\Lambda_jf\cdot\phi_y^\delta\|_q^q\dif y\right)^{1/q}\leq\|\Lambda_jf\|_q\leq2^{-\gamma j}\|f\|_{B^{\gamma}_{q,\infty}}.
\end{align*}
For the second term, since
$$
\sup_{x\in\mR^d}\int_{\mR^d}|\nabla \phi_y^\delta|^q\dif y<\infty,
$$
we have
\begin{align*}
\cR_j^2=\left(\int_{\mR^d}\|b\cdot \Lambda_ju\cdot\nabla\phi_y^\delta\|_q^q\dif y\right)^{1/q}\preceq\|b\cdot \Lambda_ju\|_q\leq\|b\|_\infty\|\Lambda_ju\|_q\leq\|b\|_{B^\beta_{p,\infty}}2^{-\gamma j}\|u\|_{B^{\gamma}_{q,\infty}}.
\end{align*}
For the third term, we have by (\ref{comm}) that, for every $\eps>0$, there exists a constant $c_\eps>0$ such that
\begin{align*}
\cR_j^3&=\left(\int_{\mR^d}\|[\Lambda_j,b\cdot\nabla]u\cdot\phi_y^\delta\|_q^q\dif y\right)^{1/q}\leq\|[\Lambda_j,b\cdot\nabla]u\|_q\\
&\preceq 2^{-\gamma j}\|b\|_{B^\beta_{p,\infty}}\|u\|_{B^{1-\beta+\gamma}_{\hat p,1}}\leq \eps2^{-\gamma j}\|u\|_{B^{\alpha+\gamma}_{q,\infty}}+c_\eps2^{-\gamma j}\|u\|_{B^{\gamma}_{q,\infty}},
\end{align*}
where in the last inequality we used the embedding theorem.
For the fourth term, we use Lemma \ref{newc}
to deduce that
\begin{align*}
\cR_j^4&=\left(\int_{\mR^d}\|[\Lambda_j,\sL^\sigma_\nu]u\cdot\phi_y^\delta\|_q^q\dif y\right)^{1/q}\leq\|[\Lambda_j,\sL^\sigma_\nu]u\|_q\\
&\preceq2^{-\gamma j}\|\sigma\|_{L^\infty_\infty(T)C^{\vartheta}_b}\|u\|_{B^{\alpha-\bar\vartheta+\gamma}_{q,\infty}}\leq \eps2^{-\gamma j}\|u\|_{B^{\alpha+\gamma}_{q,\infty}}+c_\eps2^{-\gamma j}\|u\|_{B^{\gamma}_{q,\infty}}.
\end{align*}
For the fifth term, we have by \cite[(2.19)]{C-Z} and \cite[Lemma 2.2]{C-Z} that
\begin{align*}
\cR_j^5&=\left(\int_{\mR^d}\|(\sL^\sigma_\nu-\sL^{\sigma_0}_\nu)\Lambda_ju\cdot\phi_y^\delta\|_q^q\dif y\right)^{1/q}\leq\big\|\sup_{|x-y|\leq\delta}(\sL^\sigma_\nu-\sL^{\sigma_0}_\nu)\Lambda_ju\big\|_q\\
&\preceq \delta^\vartheta\|\Delta^{\alpha/2}\Lambda_ju\|_q\preceq \delta^\vartheta2^{-\gamma j}\|u\|_{B^{\alpha+\gamma}_{q,\infty}}.
\end{align*}
Finally, to handle the last term, we write
\begin{align*}
&\sL^{\sigma_0}_\nu\Lambda_ju\cdot\phi_y^\delta-\sL^{\sigma_0}_\nu(\Lambda_ju\cdot\phi_y^\delta)\\
&=-\Lambda_ju\sL^{\sigma_0}_\nu\phi^\delta_y-\int_{\mR^d}\big(\Lambda_ju(x+z)-\Lambda_ju(x)\big)\big(\phi_y^\delta(x+z)-\phi_y^\delta(x)\big)\sigma(t,y,z)\nu(\dif z).
\end{align*}
Since
$$
\sup_{x\in\mR^d}\int_{\mR^d}\sL^{\sigma_0}_\nu\phi^\delta_y\dif y<\infty,
$$
we have
$$
\left(\int_{\mR^d}\|\Lambda_ju\sL^{\sigma_0}_\nu\phi^\delta_y\|_q^q\dif y\right)^{1/q}\preceq\|\Lambda_ju\|_q\leq 2^{-\gamma j}\|u\|_{B^\gamma_{q,\infty}}.
$$
For the second part, let
\begin{align*}
\cV:=\int_{\mR^d}\big(\Lambda_ju(x+z)-\Lambda_ju(x)\big)\big(\phi_y^\delta(x+z)-\phi_y^\delta(x)\big)\sigma(t,y,z)\nu(\dif z).
\end{align*}
Then we can deduce that
\begin{align*}
&\left(\int_{\mR^d}\!\|\cV\|_q^q\dif y\right)^{1/q}\!\!\leq\int_{\mR^d}\!\left(\int_{\mR^d}\!\|\big(\Lambda_ju(x+z)-\Lambda_ju(x)\big)\big(\phi_y^\delta(x+z)-\phi_y^\delta(x)\big)\|_q^q\dif y\right)^{1/q}\!\nu(\dif z)\\
&\leq\int_{\mR^d}\|\big(\Lambda_ju(x+z)-\Lambda_ju(x)\big)\|_q\left(\sup_{x\in\mR^d}\int_{\mR^d}|\big(\phi_y^\delta(x+z)-\phi_y^\delta(x)\big)|^q\dif y\right)^{1/q}\nu(\dif z)\\
&\preceq\|\Lambda_ju\|_q\int_{\mR^d}\left(\sup_{x\in\mR^d}\int_{\mR^d}|\big(\phi_y^\delta(x+z)-\phi_y^\delta(x)\big)|^q\dif y\right)^{1/q}\nu(\dif z)\\
&\preceq\|\Lambda_ju\|_q\leq 2^{-\gamma j}\|u\|_{B^\gamma_{q,\infty}},
\end{align*}
where in the last inequality we used the fact that
$$
\int_{\mR^d}\left(\sup_{x\in\mR^d}\int_{\mR^d}|\big(\phi_y^\delta(x+z)-\phi_y^\delta(x)\big)|^q\dif y\right)^{1/q}\nu(\dif z)<\infty.
$$
Thus, we have
$$
\cR_j^6=\left(\int_{\mR^d}\big\|\big(\sL^{\sigma_0}_\nu\Lambda_ju\cdot\phi_y^\delta -\sL^{\sigma_0}_\nu(\Lambda_ju\cdot\phi_y^\delta)\big)\big\|_q^q\dif y\right)^{1/q}\preceq2^{-\gamma j}\|u\|_{B^\gamma_{q,\infty}}.
$$
Combining the above estimates, we get (\ref{ess}). As a result, we arrive at that for some positive constant $\kappa_2>0$,
\begin{align*}
\|\Lambda_ju(t)\|_q&\leq\int_0^t-(\kappa_02^{\alpha j}+\lambda-\kappa_1)\|\Lambda_ju(s)\|_q\dif s\\
&+\int_0^t2^{-\gamma j}\kappa_2\big((\eps+\delta^\vartheta)\|u(s)\|_{B^{\alpha+\gamma}_{q,\infty}}+C_\eps\|u(s)\|_{B^\gamma_{p,\infty}}+\|f(s)\|_{B^{\gamma}_{q,\infty}}\big)\dif s.
\end{align*}
Then, choosing $\eps$ and $\delta$ small enough so that $\kappa_2(\eps+\delta^\vartheta)<1$, and following the same argument as in the proof
of Theorem \ref{drift2}.(i), we can get the desired result.

\vspace{2mm}
The conclusions (ii) and (iii) can be proved by using the same arguments as above and
following the proofs of Theorem \ref{drift2}.(ii) and Theorem \ref{drift2}.(iii), respectively, we omit the details.
\end{proof}

We also study the equation (\ref{pde}) in H\"older spaces. To this end, we will need the following lemma.

\bl\label{pos}
Assume that {\bf (H$^\nu$)} holds and for $j\geq0$, let
$$
\cA:=\Big\{u\in C_0(\mR^d): \rm{supp}\, \hat u \subseteq\{\xi: 2^{j-1}\leq|\xi|\leq 2^{j+1}\}\Big\}.
$$
Suppose that $u\in\cA$ and $|u(x_0)|=\|u\|_\infty$ for some $x_0\in\mR^d$. Then there exists a
constant $C_0$ independent of $u$ such that for every $\kappa>0$,
$$
\rm{sign}(u(x_0))\sL_{\nu}^{\kappa}u(x_0)\leq -C_02^{\alpha j}\|u\|_\infty,
$$
where $\sL_{\nu}^{\kappa}$ is defined by (\ref{ope}) with $\sigma\equiv\kappa$.
\el
\begin{proof}
Without loss of generality, we may assume that $\|u\|_\infty=1$. Then, following the same argument as in the proof of \cite[Lemma 3.4]{W-Z} (since only the positivity of $\nu$ and the change of order for the translation operator with $\sL_\nu^\kappa$ are needed), we  can show that for any $u\in C_0(\mR^d)$ with $\rm{supp}\, \hat u \subseteq\{\xi: 1/2\leq|\xi|\leq 2\}$, there exists a positive constant $c_0$ such that
$$
\rm{sign}(u(x_0))\sL_{\nu}^{\kappa}u(x_0)\leq -c_0.
$$
Now, for every $u\in\cA$ with $u(x_0)=\|u\|_\infty$, let $g(x):=u(2^{-j}x)$. It is easy to see that $\rm{supp}\, \hat g \subseteq\{\xi: 1/2\leq|\xi|\leq 2\}$ and $g(2^jx_0)=\|g\|_\infty$. Thus, we have
\begin{align*}
\sL_{\nu}^\kappa u(x_0)&=\sL_{\nu}^\kappa g(2^jx_0)\leq\int_{|z|\leq1}\big[g(2^jx_0+2^jz)-g(2^jx_0)\big]\nu_1(\dif z)\\
&\preceq 2^{\alpha j}(\sL_\nu^\kappa g)(2^jx_0)\preceq -c_02^{\alpha j},
\end{align*}
where we have used the fact that $g(2^jx_0+2^jz)-g(2^jx_0)\leq 0$. The proof is complete.
\end{proof}

The following result corresponds to the case $q=\infty$ (thus $p=\infty$) in Theorem \ref{drift3}. That is, we have the following result for equation (\ref{pde}) in H\"older spaces.

\bt\label{drift}
(i) Assume that $0<\alpha\leq1$, {\bf (H$^\sigma_1$)} holds, $b\in L^\infty\big([0,T];C_b^\beta(\mR^d)\big)$ with $\beta>1-\alpha$. Then, for any $\lambda\geq0$ and $f\in L^\infty\big([0,T];C_b^\gamma(\mR^d)\big)$ with $\gamma\in[0,\vartheta\wedge\beta]$, there exists a unique solution $u\in L^\infty\big([0,T];B^{\alpha+\gamma}_{\infty,\infty}(\mR^d)\big)$ to (\ref{pde}) with
\begin{align}
\|u\|_{L^\infty_TB^{\alpha+\gamma}_{\infty,\infty}}\leq C_1\|f\|_{L^\infty_TC_b^{\gamma}},  \label{esh1}
\end{align}
where $C_1=C(d,T,\alpha,\gamma,\|b\|_{L^\infty_TC^{\beta}_b})$ is a positive constant. Moreover, for any $\eta\in[0,\alpha+\gamma)$,
\begin{align}
\|u\|_{L^\infty_TB_{\infty,\infty}^{\eta}}\leq C_\lambda\|f\|_{L^\infty_TC_b^\gamma} , \label{essh11}
\end{align}
where $C_\lambda$ is a positive constant satisfying $C_\lambda\rightarrow0$ as $\lambda\rightarrow\infty$.

\vspace{2mm}
\noindent (ii) Assume that $\alpha\in(0,1)$, {\bf (H$^\sigma_1$)} holds, $b\in L^\infty\big([0,T];C^{1-\alpha}_b(\mR^d)\big)$ with $\|b\|_{L^\infty_TC^{1-\alpha}_b}$ small. Then, for any $\lambda\geq0$ and $f\in L^\infty\big([0,T];C^\gamma_b(\mR^d)\big)$ with $\gamma\in(0,(1-\alpha)\wedge\vartheta)$, there exists a unique solution $u\in L^\infty\big([0,T];C^{\alpha+\gamma}_b(\mR^d)\big)$ to (\ref{pde}) with
\begin{align}
\|u\|_{L^\infty_TC^{\alpha+\gamma}_b}\leq C_2\|f\|_{L^\infty_TC^\gamma_b},    \label{essh2}
\end{align}
where $C_2=C(d,T,\alpha,\gamma,\|b\|_{L^\infty_TC^{1-\alpha}_b})$ is a positive
constant. Moreover, for any $\eta\in[0,\alpha+\gamma)$, (\ref{essh11}) holds with
a positive constant $C_\lambda$ satisfying $C_\lambda\rightarrow0$ as $\lambda\rightarrow\infty$.

\vspace{2mm}
\noindent (iii) Assume that $\alpha\in(0,1)$, {\bf (H$^\sigma_1$)} holds, $b\in L^\infty\big([0,T];B^{1-\alpha}_{\infty,1}(\mR^d)\big)$ with $\|b\|_{L^\infty_TB^{1-\alpha}_{\infty,1}}$ small. Then, for any $\lambda\geq0$ and $f\in L^\infty\big([0,T];B^0_{\infty,\infty}(\mR^d)\big)$, there exists a unique solution $u\in L^\infty\big([0,T];C^{\alpha}_b(\mR^d)\big)$ to (\ref{pde}) with
\begin{align}
\|u\|_{L^\infty_TC^{\alpha}_b}\leq C_3\|f\|_{L^\infty_TB^0_{\infty,\infty}},    \label{essh6}
\end{align}
where $C_3=C(d,T,\alpha,\|b\|_{L^\infty_TC^{1-\alpha}_b})$ is a positive
constant. Moreover, for any $\eta\in[0,\alpha)$, (\ref{essh11}) holds with $\gamma=0$ and
a positive constant $C_\lambda$ satisfying $C_\lambda\rightarrow0$ as $\lambda\rightarrow\infty$.
\et

\begin{proof}
As in the proof of Theorem \ref{drift2}, we only need to prove the a priori estimates
and (\ref{essh11}). That is, we assume that $u\in L^\infty\big([0,T];C_b^{\alpha+\gamma}(\mR^d)\big)$ satisfies (\ref{pde}), we show  (\ref{esh1}), (\ref{essh11}), (\ref{essh2}) and (\ref{essh6}) hold true.

\vspace{2mm}
(i) Using the operator $\Lambda_j$ to act both sides of (\ref{pde}), we get
\begin{align*}
\p_t\Lambda_ju&=\Lambda_j(\sL^{\sigma}_\nu u)+\Lambda_j(b\cdot\nabla u)+\Lambda_jf-\lambda\Lambda_ju\\
&=\sL^{\sigma}_\nu\Lambda_ju+[\Lambda_j,\sL^\sigma_\nu]u+b\cdot\nabla\Lambda_ju+[\Lambda_j,b\cdot\nabla]u-\lambda\Lambda_ju+\Lambda_jf.
\end{align*}
For each $t\in[0,T]$,
let $x_{t,j}$ be the point such that $\Lambda_j u(t)$ reaches its maximum at $x_{t,j}$.
Without loss of generality, we may assume that
$\Lambda_j u(t, x_{t,j})>0$ (otherwise, we may consider the equation
for $-\Lambda_j u$ instead). Now, by using the fact that
$\Lambda_j u(t,x_{t,j}+y)-\Lambda_j u(t,x_{t,j})\leq0$, we get that
\begin{align*}
\sL^{\sigma}_\nu\Lambda_ju(t,x_{t,j})&=\int_{\mR^d}\big[\Lambda_j u(t,x_{t,j}+y)
-\Lambda_j u(t,x_j)\big]\sigma(t,x,z)\nu(\dif z)\\
&\leq \int_{\mR^d}\big[\Lambda_j u(t,x_{t,j}+y)
-\Lambda_j u(t,x_{t,j})\big]\kappa_0\nu(\dif z)
=\sL^{\kappa_0}_\nu\Lambda_ju(t,x_{t,j}),
\end{align*}
where $\kappa_0$ is the lower bound of $\sigma$ given in {\bf (H$^\sigma_1$)}.
According to Lemma \ref{pos}, we can further get that there exists a constant $c_0>0$ such that for all $j\geq0$,
$$
\sL^{\sigma}_\nu\Lambda_ju(t,x_{t,j})\leq -c_02^{\alpha j}\|\Lambda_ju\|_\infty.
$$
On the other hand, we can use Lemma \ref{newc} with $\bar\vartheta>0$ to get that for some $c_1>0$,
\begin{align*}
[\Lambda_j,\sL^\sigma_\nu]u(t,x_{t,j})\leq\|[\Lambda_j,\sL^\sigma_\nu]u\|_\infty\leq c_12^{-\gamma j}\|\sigma\|_{L^\infty_\infty(T)C_b^\vartheta}\|u\|_{B^{\alpha-\bar\vartheta+\gamma}_{\infty,\infty}},\quad\forall j\geq-1.
\end{align*}
For the third term, it is obvious that
$$
b\cdot\nabla\Lambda_ju(t,x_{t,j})=0.
$$
Similarly, we can use (\ref{comm}) with $p=p_1=p_2=\infty$ to deduce that for every $\gamma\in[0,\beta]$, there exists a constant $c_2>0$ such that
\begin{align*}
[\Lambda_j,b\cdot\nabla]u(t,x_{t,j})\leq\|[\Lambda_j,b\cdot\nabla]u\|_\infty\leq c_22^{-\gamma j}\|b\|_{B^\beta_{\infty,\infty}}\|u\|_{B^{1-\beta+\gamma}_{\infty,1}}.
\end{align*}
Thus, we have by \cite[Lemma 3.2]{W-Z} that for some $c_3>0$,
\begin{align*}
\frac{\dif}{\dif t}\|\Lambda_ju(t)\|_\infty&=\p_t\Lambda_ju(t,x_{t,j})\leq -c_02^{\alpha j}\|\Lambda_ju\|_\infty+c_12^{-\gamma j}\|b\|_{B^\beta_{\infty,\infty}}\|u\|_{B^{1-\beta+\gamma}_{\infty,1}}\\
&\quad+c_12^{-\gamma j}\|\sigma\|_{L^\infty_\infty(T)C_b^\vartheta}\|u\|_{B^{\alpha-\bar\vartheta+\gamma}_{\infty,\infty}}-(\lambda-c_3)\|\Lambda_ju\|_\infty+\|\Lambda_jf(t)\|_\infty.
\end{align*}
Using the same argument as in the proof of Theorem \ref{drift2}.(i).
we can get that estimates (\ref{esh1}) and (\ref{essh11}) hold true.

The proofs of (ii) and (iii) follow the same argument  as above and the procedures in the proofs of Theorem \ref{drift2}.(ii) and Theorem \ref{drift2}.(iii), respectively.
The proof is complete.
\end{proof}

Part (iii) of Theorems
\ref{drift3} and \ref{drift} will not be used later in the paper. We include these assertions since they may be of independent interest in the theory of PDEs.

In the case $\alpha=1$,
we have the following result for equation  (\ref{pde}) in Bessel potential spaces, which will be enough for us to prove the well-posedness of the martingale problem for $\sL_t$ when $\alpha=1$ and $b$ is bounded.

\bt\label{th}
Assume that $\alpha=1$, {\bf (H$^\sigma_1$)} holds, $b\in L^\infty([0,T];L^\infty(\mR^d))$ with $\|b\|_{L^\infty_\infty(T)}$ small. Then, for any $f\in L^\infty([0,T];L^q(\mR^d))$ with $2\leq q<\infty$, there exists a unique solution $u\in L^\infty([0,T];H^{1,q}(\mR^d))$ to (\ref{pde}) such that
\begin{align}
\|\p_tu\|_{L^\infty_q(T)}+\|u\|_{L^\infty_TH^{1,q}}\leq C_1\|f\|_{L^\infty_q(T)},   \label{essh3}
\end{align}
where $C_1=C_1(d,T,q,\|b\|_{L^\infty_\infty(T)})$ is a positive constant.
\et

\begin{proof}
We only prove the a priori estimate (\ref{essh3}).
By (\ref{pp}), we can use Theorem \ref{drift3}.(i) with $\alpha=1$ and $b=0$ to get that for every $f\in L^\infty([0,T];L^q(\mR^d))$ with $2\leq q<\infty$, there exists a unique solution $u\in L^\infty([0,T];B^{1+\vartheta}_{q,\infty}(\mR^d))$ to the following equation:
$$
\p_tu-\sL^\sigma_\nu u+\lambda u=f,\quad u(0,x)=0.
$$
By the embedding  (\ref{emb2}), we have that for some constant $c_0>0$,
$$
\|u\|_{L^\infty_TH^{1,q}}\leq\|u\|_{L^\infty_TB^{\alpha+\vartheta}_{q,\infty}}\leq c_0\|f\|_{L^\infty_q(T)}.
$$
It also follows from the equation itself that
$$
\|\p_tu\|_{L^\infty_q(T)}\leq c_0\|f\|_{L^\infty_q(T)}.
$$
Thus, by writing  equation (\ref{pde}) in the form
$$
\p_tu-\sL^\sigma_\nu u+\lambda u=f+b\cdot\nabla u,
$$
we can get that
\begin{align*}
\|\p_tu\|_{L^\infty_q(T)}+\|u\|_{L^\infty_TH^{1,q}}&\leq c_0\big(\|f\|_{L^\infty_q(T)}+\|b\cdot\nabla u\|_{L^\infty_q(T)}\big)\\
&\leq C_0\|f\|_{L^\infty_q(T)}+c_0\|b\|_{L^\infty_\infty(T)}\|\nabla u\|_{L^\infty_q(T)}.
\end{align*}
Consequently, if $\|b\|_{L^\infty_\infty(T)}$ is small enough so that $\|b\|_{L^\infty_\infty(T)} c_0<1$, we have
$$
\|\p_tu\|_{L^\infty_q(T)}+\|u\|_{L^\infty_TH^{1,q}}\leq c_1\|f\|_{L^\infty_q(T)}.
$$
The proof is complete.
\end{proof}

\section{Proofs of the main results}

We will first derive Krylov's estimate for solutions of SDE (\ref{sdee}) in Subsection 5.1. Then, we give the proofs of Theorems \ref{main1} and \ref{main2} in Subsection 5.2, and prove Theorem \ref{main3} in Subsection 5.3.

\subsection{Krylov's estimate}

We first give a generalization of It\^o's formula for SDE (\ref{sdee}), which will be used several times below.

\bl\label{ito}
Let $X_t$ solves (\ref{sdee}) and $f\in L^\infty(\mR_+;C^{\gamma}_b(\mR^d))$ with $\gamma=1$ when $\alpha\in(0,1)$ and $\gamma>1$ when $\alpha=1$, and $\p_tf\in L^\infty(\mR_+;C_b(\mR^d))$. Then we have for every $t\geq 0$,
\begin{align*}
f(t,X_t)-f(0,x)-\int_0^t\!(\p_s+\sL_s) f(s,X_s)\dif s=&\int_0^t\!\!\!\int_0^{\infty}\!\!\!\!\int_{\mR^d}\Big(f\big(s,X_{s-}+1_{[0,\sigma(s,X_{s-},z)]}(r)z\big)\\
&\quad-f(s,X_s)\Big)\widetilde\cN(\dif z\times\dif r\times \dif s).
\end{align*}
\el
\begin{proof}
Let $\rho\in C^{\infty}_0(\mR_+\times\mR^d)$ be such that $\int_{\mR_+\times\mR^d}\rho(t,x)\dif x\dif t=1$. Define $\rho_n(t,x):=n^{d+1}\rho(nt,nx)$ and
\begin{align}
f_n(t,x):=\int_{\mR_+\times\mR^d}f(s,y)\rho_n(t-s,x-y)\dif y\dif s.\label{mo}
\end{align}
Hence, we have $f_n\in C_b(\mR_+;C^2_b(\mR^d))\cap C^1_b(\mR_+;C_b(\mR^d))$ with $\|f_n(t)\|_{C^\gamma_b}\leq \|f(t)\|_{C^\gamma_b}$, $\|\p_tf_n\|_\infty\leq\|\p_tf\|_\infty$ and $\|f_n-f\|_{C^{\gamma'}_b}\rightarrow0$ for every $\gamma'<\gamma$. By using It\^o's formula for $f_n(t,X_t)$, we can get
\begin{align*}
f_n(t,X_t)-f_n(0,x)-\int_0^t\!(\p_s+\sL_s) f_n(s,X_s)\dif s=&\int_0^t\!\!\!\int_0^{\infty}\!\!\!\!\int_{\mR^d}\big[f_n\big(s,X_{s-}+1_{[0,\sigma(X_{s-},z)]}(r)z\big)\\
&\quad-f_n(s,X_{s-})\big]\tilde \cN(\dif z\times\dif r\times \dif s).
\end{align*}
Now we are going to let $n\to\infty$ on the both sides of the above equality. It is easy to see that for every $\omega$ and $x\in\mR^d$,
$$
f_n(t,X_t)-f_n(0,x)\rightarrow f(t,X_t)-f(0,x),\quad \text{as}\,\,n\rightarrow\infty.
$$
Since
\begin{align*}
|f_n(t,x+z)-f_n(t,x)-z\cdot\nabla f_n(t,x)|\leq c|z|^{\gamma}\|f_n\|_{C^\gamma_b}\leq c|z|^{\gamma}\|f\|_{C^\gamma_b},
\end{align*}
we can get by the dominated convergence theorem that for every $\omega$,
$$
\int_0^t\!(\p_s+\sL_s) f_n(s,X_s)\dif s\rightarrow\int_0^t\!(\p_s+\sL_s) f(s,X_s)\dif s,\quad \text{as}\,\,\,n\rightarrow\infty.
$$
Finally, by the isometry formula, we have
\begin{align*}
&\mE\bigg|\!\int_0^t\!\!\!\int_0^{\infty}\!\!\!\!\int_{\mR^d}\Big[f_n\big(s,X_{s-}+1_{[0,\sigma(X_{s-},z)]}(r)z\big)-f_n(s,X_s)\\ &\qquad-f\big(s,X_{s-}+1_{[0,\sigma(X_{s-},z)]}(r)z\big)+f(s,X_{s-})\Big]\tilde \cN(\dif z\times\dif r\times \dif s)\bigg|^2\\
&=\mE\int_0^t\!\!\!\int_{\mR^d}\!\!\int_0^{\infty}1_{[0,\sigma(X_s,z)]}(r)\big|f_n(s,X_s+z)-f_n(s,X_s)\\
&\qquad\qquad\qquad\qquad\qquad-f(s,X_s+z)+f(s,X_s)\big|^2\dif r\nu(\dif z)\dif s\\
&\leq c\!\int_0^t\!\!\!\int_{\mR^d}\mE\big|f_n(s,X_s+z)-f_n(s,X_s)\\
&\qquad\qquad\qquad\qquad\qquad-f(s,X_s+z)+f(s,X_s)\big|^2\nu(\dif z)\dif s\rightarrow0,\quad \text{as}\,\,n\rightarrow\infty,
\end{align*}
where in the last step we have used the fact that $\sigma$ is bounded, $\|f_n\|_{C^\gamma_b}\leq \|f\|_{C^\gamma_b}$ and the dominated convergence theorem again. The proof is complete.
\end{proof}

Now, we prove the following Krylov type estimate for solutions of SDE (\ref{sdee}), which will play an important role in proving  the pathwise uniqueness of strong solutions.

\bl\label{kkr}
Assume that $0<\alpha\leq1$, {\bf (H$^\sigma_1$)} holds true, $b\in L^\infty\big([0,T];B^\beta_{p,\infty}(\mR^d)\big)$
with $\beta\in (1-\alpha, 1]$ and $d/(\alpha+\beta-1)\vee2<p\leq\infty$.
Let $X_t(x)$ solves SDE (\ref{sdee}). Then, for every $f\in L^\infty\big([0,T];B^0_{q,\infty}(\mR^d)\big)$
with $q\in (\frac{d}\alpha\vee2, p]$, we have
\begin{align}
\sup_{x\in\mR^d}\mE\left(\int_0^T\!f\big(s,X_s(x)\big)\dif s\right)\leq C_d\|f\|_{L^\infty_TB^0_{q,\infty}},    \label{kry}
\end{align}
where $C_d>0$ is a constant independent of $x$.
\el
\begin{proof}
By a standard density argument, it suffices to prove the lemma for $f\in C_0^\infty(\mR_+\times\mR^d)$. Without loss of generality, we may assume $\beta$ is slightly bigger than $1-\alpha$ so that $\vartheta>\beta+\alpha-1$.
Without loss of generality, we may also assume that $\vartheta\leq \beta$.
Let $\rho\in C_0^\infty(\mR^d)$ be such that $\int_{\mR^d}\rho(x)dx=1$. Define
$\rho_n(x):=n^d\rho(nx)$ and
\begin{align*}
\sigma_n(t,x,z):=\int_{\mR^d}\sigma(t,y,z)\rho_n(x-y)\dif y.
\end{align*}
Then, $\sigma_n\in L^\infty\big([0,T]\times\mR^d_z;C_b^\infty(\mR^d_x)\big)$.
For fixed $T>0$, let $u_n$ be the solution to the following backward equation: for $t\in[0,T]$,
\begin{align}
\p_tu_n+\sL^{\sigma_n}_\nu u_n+b\cdot\nabla u_n-\lambda u_n+f=0,\quad u_n(T,x)=0.  \label{tem}
\end{align}
Thus, according to Theorem \ref{drift3}.(i), we have
for any $0\leq\gamma\leq\beta$,
\begin{align}
\|u_n\|_{L^\infty_TB^{\alpha+\gamma}_{p,\infty}}\leq c(n)\|f\|_{L^\infty_TB^\gamma_{p,\infty}}, \quad\forall n\geq1,  \label{gaa}
\end{align}
where $c(n)$ is a positive constant depending on $\|b\|_{L^\infty_TB^\beta_{p,\infty}}$ and
$\|\sigma_n\|_{L^\infty_\infty(T)C_b^1}$.

Taking $\gamma=\beta$ in \eqref{gaa} and noticing that
$$
B^{\alpha+\beta}_{p,\infty}\subseteq C_b^\eta(\mR^d)\quad\text{with}\quad\eta>1,
$$
we get
$$
u_n\in L^{\infty}([0,T];C_b^\eta(\mR^d)).
$$
By the equation (\ref{tem}) itself, we also have
$$
\p_tu_n=-\sL^{\sigma_n}_\nu u_n-b\cdot\nabla u_n+\lambda u_n-f\in L^{\infty}([0,T];C_b(\mR^d)).
$$
Consequently, by Lemma \ref{ito}, we can use It\^o's formula and take expectation to get that
\begin{align*}
\mE u_n(t,X_t)&=u_n(0,x)+\mE\left(\int_0^t\big(\p_su_n+\sL^{\sigma_n}_\nu u_n+b\cdot\nabla u_n\big)(s,X_s)\dif s\right)\\
&=u_n(0,x)+\mE\left(\int_0^t\lambda u_n(s,X_s)\dif s\right)-\mE\left(\int_0^tf(s,X_s)\dif s\right)\\
&+\mE\left(\int_0^t\!\!\!\int_{\mR^d}[u_n(s,X_s+z)-u_n(s,X_s)]\big(\sigma_n(s,X_s,z)-\sigma(s,X_s,z)\big)\nu(\dif z)\dif s\right).
\end{align*}
This in turn yields that for $\hat\eta=\alpha+\vartheta-\frac{d}{p}$,
\begin{align*}
\sup_{x\in\mR^d}\mE\left(\int_0^tf\big(s,X_s(x)\big)\dif s\right)&\leq c_{\lambda,t}\|u_n\|_\infty\\
&\quad+\mE\left(\int_0^t|\sigma_n(s,X_s,z)-\sigma(s,X_s,z)|\dif s\right)\|u_n\|_{L^\infty_TC_b^{\hat\eta}}.
\end{align*}
Recall that $q>d/\alpha$ by assumption.
It follows from Theorem \ref{drift3}.(i) that for $0\le \gamma\le \vartheta$ and
$\tilde q\in (\frac{d}\alpha\vee2, p]$,
\begin{align}\label{gaa2}
\|u_n\|_{L^\infty_TB^{\alpha+\gamma}_{\tilde q,\infty}}\leq c_1\|f\|_{L^\infty_TB^\gamma_{\tilde q,\infty}}, \quad\forall n\geq1,
\end{align}
where $c_1>0$ is a constant depending on $\|b\|_{L^\infty_TB^\beta_{p,\infty}}$ and
$\|\sigma\|_{L^\infty_\infty(T)C_b^\vartheta}$.
Taking $\gamma=0$ and $\tilde q=q$ in (\ref{gaa2})
and applying Sobolev's embedding theorem (\ref{embb}), we get
$$
\|u_n\|_\infty\leq \|u_n\|_{L^\infty_TB^{\alpha}_{q,\infty}}\leq c_2\|f\|_{L^\infty_TB^{0}_{q,\infty}}, \quad\forall n\geq1,
$$
where $c_2$ depends only on $\|b\|_{L^\infty_TB^\beta_{p,\infty}}$ and $\|\sigma\|_{L^\infty_TC_b^\vartheta}$.
On the other hand, we can take $\gamma=\vartheta$
and $\tilde q=p>d/(\alpha+\beta-1)>d/\vartheta$ in (\ref{gaa2}), and using Sobolev's embedding theorem (\ref{embb}) again  to get that
$$
\|u_n\|_{L^\infty_TC_b^{\hat\eta}}\leq \|u_n\|_{L^\infty_TB^{\alpha+\vartheta}_{p,\infty}}\leq c_3\|f\|_{L^\infty_TB^{\vartheta}_{p,\infty}}, \quad\forall n\geq1,
$$
where $c_3$ also depends only on $\|b\|_{L^\infty_TB^\beta_{p,\infty}}$ and $\|\sigma\|_{L^\infty_TC_b^\vartheta}$ with $\vartheta>0$.
Thus, by the dominated convergence theorem,
we can let $n\rightarrow\infty$ to get that for every $t\in[0,T]$,
\begin{align*}
&\sup_{x\in\mR^d}\mE\left(\int_0^tf\big(s,X_s(x)\big)\dif s\right)\leq c_{\lambda,t}c_2\|f\|_{L^\infty_TB^{0}_{q,\infty}}\\
&\qquad\qquad+c_3\lim_{n\rightarrow\infty}\mE\left(\int_0^t|\sigma_n-\sigma|(s,X_s,z)\dif s\right)\|f\|_{L^\infty_TB^{\vartheta}_{p,\infty}}=c_{\lambda,t}c_2\|f\|_{L^\infty_TB^{0}_{q,\infty}}.
\end{align*}
The proof is complete.
\end{proof}

In the case that $\alpha=1$, we also need to derive Krylov's estimate for $X_t$ when
$b\in L^\infty(\mR_+\times\mR^d)$. To this end, we consider the following quasi-linear backward non-local parabolic equation:
\begin{align}
\p_tu+\sL^\sigma_\nu u+\kappa|\nabla u|-\lambda u+f=0, \quad u(T,x)=0,   \label{qua}
\end{align}
where $\kappa>0$ is a constant. We prove the following result.

\bl
Assume $\alpha=1$ and {\bf (H$^\sigma_1$)} holds. Then, for every $f\in L^\infty([0,T];C_0^\infty(\mR^d))$ and $\kappa>0$ small, there exists a unique solution $u\in L^\infty([0,T];C_b^{1+\vartheta}(\mR^d))$ to equation (\ref{qua}) such that for any $0<\gamma\leq\vartheta$,
\begin{align}
\|u\|_{L^\infty_TC_b^{1+\gamma}}\leq C_0\|f\|_{L^\infty_TC_b^{\gamma}},  \label{pri}
\end{align}where $C_0$ is a positive constant.
Moreover, for any $p>d$, we also have
\begin{align}
\|u\|_\infty\leq C_0\|f\|_{L^\infty_p(T)}.  \label{kk}
\end{align}

\el
\begin{proof}
We only need to prove the a priori estimate (\ref{pri}). In fact, by using Theorem \ref{drift}.(i) with $\alpha=1$ and $b=0$, we have that for every $f\in L^\infty([0,T];C_0^\infty(\mR^d))$, there exists a unique solution $u\in L^\infty([0,T];C_b^{1+\vartheta}(\mR^d))$ to the following equation:
$$
\p_tu+\sL^\sigma_\nu u-\lambda u+f=0, \quad u(T,x)=0.
$$
Thus for any solution $u$ to (\ref{qua}), we have that for any $\gamma\in(0,\vartheta]$, there exists a constant $c_0>0$ such that
$$
\|u\|_{L^\infty_TC_b^{1+\gamma}}\leq c_0\big(\|f\|_{L^\infty_TC_b^{\gamma}}+\kappa\||\nabla u|\|_{L^\infty_TC_b^{\gamma}}\big).
$$
Notice that for $u\in C_b^{1+\gamma}(\mR^d)$, we have
$$
|\nabla u|\in C_b^{\gamma}(\mR^d).
$$
Consequently, we have
$$
\|u\|_{L^\infty_TC_b^{1+\gamma}}\leq c_0\|f\|_{L^\infty_TC_b^{\gamma}}+c_0\kappa\|u\|_{L^\infty_TC_b^{1+\gamma}},
$$
which implies the desired result by choosing $\kappa$ small enough so that $c_0\kappa<1$.

On the other hand, by repeating the arguments in the proof of Theorem \ref{th}, we can also show that $u\in L^\infty([0,T];H^{1,p}(\mR^d))$ for any $p\geq2$ and
$$
\|u\|_{L^\infty_TH^{1,p}}\leq c_1\|f\|_{L^\infty_p(T)},
$$
where $c_1>0$ is a constant. As a result of Sobolev's embedding, we have for any $p>d$,
$$
\|u\|_\infty\leq \|u\|_{L^\infty_TH^{1,p}}\preceq\|f\|_{L^\infty_p(T)}.
$$
The proof is complete.
\end{proof}

With the above result in hand, we show the following Krylov's estimate for $X_t$ in the particular case of $\alpha=1$.

\bl
Assume that $\alpha=1$, {\bf (H$^\sigma_1$)} holds and $b\in L^\infty(\mR_+\times\mR^d)$ with $\|b\|_\infty$ small. Let $X_t(x)$ solves SDE (\ref{sdee}). Then for any $f\in L^\infty([0,T];L^p(\mR^d))$ with $p>d$,
\begin{align}
\sup_{x\in\mR^d}\mE\left(\int_0^tf\big(s,X_s(x)\big)\dif s\right)\leq C_d\|f\|_{L^\infty_p(T)},    \label{kkk}
\end{align}
where $C_d>0$ is a constant independent of $x$.
\el
\begin{proof}
It suffices to prove (\ref{kkk}) for all $f\in C_0^\infty(\mR_+\times\mR^d)$. Let $u$ be the solution to equation (\ref{qua})
with $\kappa\geq\|b\|_\infty$.
Then we have $u\in L^\infty([0,T];C_b^{1+\vartheta}(\mR^d))$ and $\p_tu\in L^\infty([0,T];C_b(\mR^d))$. Thus, by Lemma \ref{ito}, we can use It\^o's formula and take expectation to get that
\begin{align*}
\mE u(t,X_t)&=u(0,x)+\mE\left(\int_0^t\big(\p_su+\sL^{\sigma}_\nu u+b\cdot\nabla u\big)(s,X_s)\dif s\right)\\
&\leq u(0,x)+\mE\left(\int_0^t\big(\p_su+\sL^{\sigma}_\nu u+\kappa\cdot|\nabla u|\big)(s,X_s)\dif s\right)\\
&=u(0,x)-\mE\left(\int_0^tf(s,X_s)\dif s\right).
\end{align*}
This in turn yields by (\ref{kk}) that
\begin{align*}
\sup_{x\in\mR^d}\mE\left(\int_0^tf\big(s,X_s(x)\big)\dif s\right)\leq 2\|u\|_\infty\leq C_0\|f\|_{L^\infty_p(T)}.
\end{align*}
The proof is complete.
\end{proof}

\subsection{Proofs of Theorems \ref{main1} and  \ref{main2}}

Under our assumptions, the existence of the martingale solutions is known, see \cite {Ba} or \cite[Theorem 4.1]{Kurz2}. Moreover, by \cite[Theorem 2.3]{Kurz2}, we know that the martingale solution for the operator $\sL_t$ is equivalent to the weak solution for SDE (\ref{sdee}). Thus, we shall focus on proving the uniqueness of weak solutions for SDE (\ref{sdee}).
Let us first give:
\begin{proof}[Proof of Theorem \ref{main1}]
Let $f\in C_0^\infty(\mR_+\times\mR^d)$ and $u$
be the solution to the following backward equation:
\begin{align}
\p_tu+\sL^{\sigma}_\nu u+b\cdot\nabla u+f=0,\quad u(T,x)=0.  \label{equ}
\end{align}
By the assumption that  $b\in L^\infty\big([0,T];B^\beta_{p,\infty}(\mR^d)\big)$ with $\beta>1-\alpha$ and $p>d/(\alpha+\beta-1)$, we have by Theorem \ref{drift3}.(i) that $u\in L^\infty\big([0,T];B^{\alpha+\beta}_{p,\infty}(\mR^d)\big)$. In particular, we have by (\ref{embb}) that
\begin{align*}
u\in L^\infty\big([0,T];C_b^\eta(\mR^d)\big)\quad\text{with}\quad\eta>1\quad\text{and}\quad\p_tu\in L^\infty\big([0,T];C_b(\mR^d)\big).
\end{align*}
According to Lemma \ref{ito}, for any weak solution $X_t$ to SDE (\ref{sdee}), we can use Ito's formula and take expectation to get that
$$
-u(0,x)=\mE\left(\int_0^T\!f(s,X_s)\dif s\right).
$$
Note that the left hand side is independent of $X_t$. This in particular implies the uniqueness of the weak solution.
\end{proof}

The key point in the above proof of Theorem \ref{main1} is that the solution $u$ to (\ref{equ}) is regular enough for us to aplly the It\^o's formula. However, such an argument can not be easily modified to prove Theorem \ref{main2}, since under the conditions in Theorem \ref{main2}, we can not get a control on $\nabla u$ with proper norms. We divide the proof of Theorem \ref{main2} into two cases.

\begin{proof}[Proof of Theorem \ref{main2} in the case $\alpha\in(0,1)$]
Let $u$ be the solution to (\ref{equ}). Then, by Theorem \ref{drift}.(ii), we have $u\in L^\infty\big([0,T]; C_b^{\alpha+\gamma}(\mR^d)\big)$ with any $0<\gamma<\vartheta$.
Using the operator $\Lambda_j$ to act on both sides of the equation, we can get
\begin{align*}
\p_t\Lambda_ju+\sL^{\sigma}_\nu \Lambda_ju+b\cdot\nabla \Lambda_ju+\Lambda_jf+[\Lambda_j,\sL^\sigma_\nu]u+[\Lambda_j,b]\nabla u=0.
\end{align*}
Applying It\^o's formula for $\Lambda_ju$ and taking expectation,  we have
\begin{align*}
-\Lambda_ju(0,x)=\mE\left(\int_0^T\Lambda_jf(s,X_s)\dif s\right)+\mE\left(\int_0^T\Big([\Lambda_j,\sL^\sigma_\nu]u+[\Lambda_j,b]\nabla u\Big)(s,X_s)\dif s\right).
\end{align*}
Summing over $j$ and using (\ref{id}), we can get
\begin{align*}
-u(0,x)=\mE\left(\int_0^Tf(s,X_s)\dif s\right)+\sum_{j\geq-1}\mE\left(\int_0^T\Big([\Lambda_j,b]\nabla u+[\Lambda_j,\sL^\sigma_\nu]u\Big)(s,X_s)\dif s\right).
\end{align*}
We now show that the second term on the right hand side equals zero. In fact, since
$$
\sum_{j\geq-1}[\Lambda_j,f]g=\sum_{j\geq-1}\big(\Lambda_j(fg)-f\Lambda_j g\big)=fg-fg=0,
$$
we only need to ensure that we can change the order of summation and integration. Noticing that  by (\ref{comm}) and Lemma \ref{newc}, we have
$$
\sum_{j\geq-1}\|[\Lambda_j,b]\nabla u\|_\infty\preceq2^{-\gamma j}\|b\|_{C_b^{1-\alpha}}\|u\|_{B^{\alpha+\gamma}_{\infty,\infty}}<\infty
$$
and
$$
\sum_{j\geq-1}\|[\Lambda_j,\sL^\sigma_\nu] u\|_\infty\preceq2^{-(\vartheta-\bar\vartheta+\gamma)j}\|\sigma\|_{L^\infty_\infty C_b^{\vartheta}}\|u\|_{B^{\alpha-\bar\vartheta+\gamma}_{\infty,\infty}}<\infty.
$$
The desired result follows immediately.
\end{proof}

\begin{proof}[Proof of Theorem \ref{main2} in the case $\alpha=1$]
For every $f\in C_0^\infty(\mR_+\times\mR^d)$, let $u$ be the solution to equation (\ref{equ}).
Let $u_n$ be the mollification of $u$ defined the same way as in (\ref{mo}) and define
$$
f_n:=\p_tu_n+\sL^\sigma_\nu u_n+b\cdot\nabla u_n.
$$
Then, by (\ref{essh3}), we have that as $n\rightarrow\infty$,
$$
\|f_n-f\|_{L^\infty_p(T)}\leq C_0\|\p_t(u_n-u)+\sL^\sigma_\nu(u_n-u)+b\cdot\nabla(u_n-u)\|_{L^\infty_p(T)}\rightarrow0.
$$
Using It\^o's formula for $u_n(t,X_t)$ and taking expectation, we  get
$$
-u_n(0,x)=\mE\left(\int_0^Tf_n(s,X_s)\dif s\right).
$$
Letting $n\rightarrow\infty$ on both sides of the above equality, it is easy to see that
$u_n(0,x)\rightarrow u(0,x)$. On the other hand, applying Krylov's estimate (\ref{kkk}),  we get
$$
\mE\left(\int_0^T(f_n-f)(s,X_s)\dif s\right)\leq C_0\|f_n-f\|_{L^\infty_p(T)}\rightarrow0.
$$
Thus we have
$$
-u(0,x)=\mE\left(\int_0^T\!f(s,X_s)\dif s\right).
$$
The proof is complete.
\end{proof}

\subsection{Proof of Theorem \ref{main3}}
Throughout this subsection, we assume that
$b$ satisfies {\bf (H$^b_3$)}.
Now, let $u$ be the solution to the following backward equation: for $t\in[0,T]$,
$$
\p_tu(t,x)+\sL^{\sigma}_\nu u(t,x)+b(t,x)\cdot\nabla u(t,x)-\lambda u(t,x)+b(t,x)=0,\quad u(T,x)=0,
$$
where $\lambda>0$ is a constant. Note that by  assumption,
$$
p>2d/\alpha>d/(\alpha+\beta-1).
$$
Thus, we have by Theorem \ref{drift3}.(i) that
$$
u\in L^\infty\big([0,T];B^{\alpha+\beta}_{p,\infty}(\mR^d)\big).
$$
Moreover, in view of (\ref{embb}) and (\ref{ess11}), we can take $\lambda$ large enough so that
\begin{align*}
\|u\|_\infty+\|\nabla u\|_\infty\leq1/2.
\end{align*}
Define
\begin{align}
\Phi_t(x):=x+u(t,x). \label{phi}
\end{align}
It is easy to see that
$$
\frac{1}{2}|x-y|\leq\big|\Phi_t(x)-\Phi_t(y)\big|\leq \frac{3}{2}|x-y|,
$$
which implies that for each $t\in[0,T]$, the map $x\rightarrow\Phi_t(x)$ is a $C^1$-diffeomorphism and
\begin{align}
\frac{1}{2}\leq \|\nabla\Phi\|_{\infty},\|\nabla\Phi^{-1}\|_{\infty}\leq 2,   \label{upd}
\end{align}
where we use $\Phi^{-1}_t(\cdot)$ to denote the inverse function of $\Phi_t(\cdot)$.
The following result is known as Zvonkin's transform, which transform the original SDE into a new one with better coefficients.

\bl\label{zvo}
Let $\Phi_t(x)$ be defined by (\ref{phi}) and $X_t$ solve SDE (\ref{sdee}). Then $Y_t:=\Phi_t(X_t)$ satisfies the following SDE:
\begin{align}
Y_t&=\Phi_0(x)+\int_0^t\!\tilde b(s,Y_s)\dif s+\int_0^t\!\!\!\int_0^{\infty}\!\!\!\!\int_{|z|\leq 1} g_s(Y_{s-},z)1_{[0,\tilde\sigma(s,Y_{s-},z)]}(r)\widetilde \cN(\dif z\times\dif r\times \dif s)\no\\
&\quad+\int_0^t\!\!\!\int_0^{\infty}\!\!\!\!\int_{|z|> 1}g_s(Y_{s-},z)1_{[0,\tilde\sigma(s,Y_{s-},z)]}(r)\cN(\dif z\times\dif r\times \dif s), \label{sde3}
\end{align}
where
\begin{align}
\tilde b(t,x):=\lambda u\big(t,\Phi^{-1}_t(x)\big)-\int_{|z|>1}\!\Big[u\big(t,\Phi^{-1}_t(x)+z\big)-u\big(t,\Phi^{-1}_t(x)\big)\Big]\tilde\sigma(t,x,z)\nu(\dif z)  \label{newb}
\end{align}
and
\begin{align}
g_t(x,z):=\Phi_t\big(\Phi^{-1}_t(x)+z\big)-x,\quad \tilde\sigma(t,x,z):=\sigma\big(t,\Phi^{-1}_t(x),z\big). \label{newg}
\end{align}
\el
\begin{proof}
The conclusion follows by applying Lemma \ref{ito} and It\^o's formula to $X_t+u(t,X_t)$.  We omit the details.
\end{proof}

Now we are in a position to give

\begin{proof}[Proof of Theorem \ref{main3}]
Let $X^1_t$ and $X^2_t$ be two strong solutions for SDE (\ref{sdee}) both starting from $x\in\mR^d$, and set
$$
Y^1_t:=\Phi_t(X^1_t),\quad Y^2_t:=\Phi_t(X^2_t).
$$
Then, by Lemma \ref{zvo}, each $Y_t^i$ satisfies (\ref{sde3}).
Since the map $\Phi$ is one-to-one and by the classical interlacing technique, we only need to prove that the following SDE has a unique strong solution:
\begin{align*}
\dif Y_t&=\tilde b(t,Y_t)\dif t+\int_0^{\infty}\!\!\!\int_{|z|\leq 1} g_t(Y_{t-},z)1_{[0,\tilde\sigma(t,Y_{t-},z)]}(r)\widetilde \cN(\dif z\times\dif r\times \dif t),
\end{align*}
where $\tilde b$ and $g$ are defined by (\ref{newb}) and (\ref{newg}), respectively. Let $Z_t:=Y_t^1-Y_t^2$, then $Z_t$ satisfies the following equation:
\begin{align*}
Z_t
&=\int_0^{t}\Big[\tilde b(s,Y^1_s)-\tilde b(s,Y^2_s)\Big]\dif s+\int_0^{t}\!\!\!\int_0^{\infty}\!\!\!\!\int_{|z|\leq 1}\Big[g_s(Y^1_{s-},z)1_{[0,\tilde\sigma(s,Y^1_{s-},z)]}(r)\\
&\quad\qquad\quad\qquad-g_s(Y^2_{s-}, z)1_{[0,\tilde\sigma(s,Y^2_{s-},z)]}(r)\Big]\widetilde \cN(\dif z\times\dif r\times \dif s)=:\cJ_{t}^1+\cJ_{t}^2.
\end{align*}
Recall that we assume $\sigma$ satisfies (\ref{a1}). Set
$$
A_1(t):=\int_0^t\Big(1+\varrho(X^1_s)+\varrho(X^2_s)\Big)
\dif s.
$$
By (\ref{upd}), we have that for almost all $\omega$ and every stopping time $\eta$,
\begin{align*}
\sup_{t\in[0,\eta]}\left|\cJ_{t}^1\right|\leq c_1\!\int_0^{\eta}|Z_s|\cdot\Big(1+\varrho(X^1_s)+\varrho(X^2_s)\Big)\dif s=c_1\!\int_0^{\eta}|Z_s| \dif A_1(s).
\end{align*}
For the second term, write
\begin{align*}
\cJ^2_{t}&=\int_0^{t}\!\!\!\int_0^{\infty}\!\!\!\!\int_{|z|\leq 1}1_{[0,\tilde\sigma(s,Y^1_{s-},z)\wedge\tilde\sigma(s,Y^2_{s-},z)]}(r)\Big[g_s(Y^1_{s-},z)1_{[0,\tilde\sigma(s,Y^1_{s-},z)]}(r)\\
&\qquad\qquad\qquad\qquad\qquad\qquad\qquad-g_s(Y^2_{s-}, z)1_{[0,\tilde\sigma(s,Y^2_{s-},z)]}(r)\Big]\widetilde \cN(\dif z\times\dif r\times \dif s)\\
&\quad+\int_0^{t}\!\!\!\int_0^{\infty}\!\!\!\!\int_{|z|\leq 1}1_{[\tilde\sigma(s,Y^1_{s-},z)\vee\tilde\sigma(s,Y^2_{s-},z),\infty]}(r)\Big[g_s(Y^1_{s-},z)1_{[0,\tilde\sigma(s,Y^1_{s-},z)]}(r)\\
&\qquad\qquad\qquad\qquad\qquad\qquad\qquad-g_s(Y^2_{s-}, z)1_{[0,\tilde\sigma(s,Y^2_{s-},z)]}(r)\Big]\widetilde \cN(\dif z\times\dif r\times \dif s)\\
&\quad+\int_0^{t}\!\!\!\int_0^{\infty}\!\!\!\!\int_{|z|\leq 1}1_{[\tilde\sigma(s,Y^1_{s-},z)\wedge\tilde\sigma(s,Y^2_{s-},z),\tilde\sigma(s,Y^1_{s-},z)\vee\tilde\sigma(s,Y^2_{s-},z)]}(r) \Big[g_s(Y^1_{s-},z)1_{[0,\tilde\sigma(s,Y^1_{s-},z)]}(r)\\
&\qquad\qquad\qquad\qquad\qquad\qquad\qquad-g_s(Y^2_{s-}, z)1_{[0,\tilde\sigma(s,Y^2_{s-},z)]}(r)\Big]\widetilde \cN(\dif z\times\dif r\times \dif s)\\
&=:\cJ^{21}_{t}+\cJ^{22}_{t}+\cJ^{23}_{t}.
\end{align*}
We proceed to estimate each of the three terms above. First, for $\cJ^{21}_{t}$, we use Doob's $L^2$-maximal inequality to deduce that for any stopping time $\eta$,
\begin{align*}
&\mE\left[\sup_{t\in[0,\eta]}|\cJ^{21}_{t}|\right]\\
&\leq \mE\Bigg(\int_0^{\eta}\!\!\!\!\int_0^{\infty}\!\!\!\!\int_{|z|\leq 1}1_{[0,\tilde\sigma(s,Y^1_s,z)\wedge\tilde\sigma(s,Y^2_s,z)]}(r)\big|g_s(Y^1_s,z)-g_s(Y^2_s, z)\big|^2\dif r\nu(\dif z) \dif s\Bigg)^{\frac{1}{2}}\\
&=\mE\Bigg(\int_0^{\eta}\!\!\!\!\int_{|z|\leq 1}\big[\tilde\sigma(s,Y^1_s,z)\wedge\tilde\sigma(s,Y^2_s,z)\big] \cdot\big|g_s(Y^1_s,z)-g_s(Y^2_s, z)\big|^2\nu(\dif z) \dif s\Bigg)^{\frac{1}{2}}.
\end{align*}
Note that by (\ref{sh}),
$$
g_s(Y^1_s,z)=\sT_z\Phi_s(X^1_s).
$$
Thus, if we set
$$
A_2(t):=\int_0^t\!\!\!\int_{|z|\leq 1}\!\Big(\cM|\nabla\cJ_zu|(s,X^1_s)+\cM|\nabla\cJ_zu|(s,X^2_s)\Big)^2\nu(\dif z) \dif s,
$$
we then have by (\ref{w11}) and the fact that $\tilde\sigma$ is bounded  that
\begin{align*}
\mE\left[\sup_{t\in[0,\eta]}|\cJ^{21}_{t}|\right]&\leq c_2\mE\Bigg(\int_0^{\eta}|Z_s|^2\!\int_{|z|\leq 1}\!\Big(\cM|\nabla\cJ_zu|(X^1_s)+\cM|\nabla\cJ_zu|(X^2_s)\Big)^2\nu(\dif z) \dif s\Bigg)^{\frac{1}{2}}\\
&=c_2\mE\Bigg(\int_0^{\eta}|Z_s|^2\dif A_2(s)\Bigg)^{\frac{1}{2}}.
\end{align*}
Next, it is easy to see that for any $t\geq 0$,
$$
\cJ^{22}_{t}\equiv0.
$$
Finally, we use the $L^1$-estimate (see \cite[p.174]{Kurz3} or \cite[p.157]{Kurz2}) to control the third term:
\begin{align*}
\mE\left[\sup_{t\in[0,\eta]}|\cJ^{23}_{t}|\right]&\leq 2\mE\!\int_0^{\eta}\!\!\!\!\int_0^{\infty}\!\!\!\!\int_{|z|\leq 1}1_{[\tilde\sigma(s,Y^1_s,z)\wedge\tilde\sigma(s,Y^2_s,z),\tilde\sigma(s,Y^1_s,z)\vee\tilde\sigma(s,Y^2_s,z)]}(r)\\
&\qquad\qquad\times\big|g_s(Y^1_s,z)1_{[0,\tilde\sigma(s,Y^1_s,z)]}(r)-g_s(Y^2_s, z)1_{[0,\tilde\sigma(s,Y^2_s,z)]}(r)\big|\nu(\dif z)\dif r \dif s\\
&\leq2\mE\!\int_0^{\eta}\!\!\!\int_{|z|\leq 1}|\tilde\sigma(s,Y^1_s,z)-\tilde\sigma(s,Y^2_s,z)|\\
&\qquad\qquad\qquad\qquad\times\Big(|g_s(Y^1_s,z)|+|g_s(Y^2_s,z)|\Big)\nu(\dif z)\dif s.
\end{align*}
Since
$$
|g_s(x,z)|=\big|\Phi_s\big(\Phi^{-1}_s(x)+z\big)-\Phi_s\big(\Phi^{-1}_s(x)\big)\big|\leq \frac{3}{2}|z|,
$$
taking into account of (\ref{a1}), we get
\begin{align*}
\mE\left[\sup_{t\in[0,\eta]}|\cJ^{23}_{t}|\right]&\leq c_3\mE\!\int_0^{\eta}\!\!\!\int_{|z|\leq 1}|\tilde\sigma(s,Y^1_s,z)-\tilde\sigma(s,Y^2_s,z)|\cdot|z|\nu(\dif z)\dif s\\
&\leq c_3\mE\!\int_0^{\eta}|Z_s|\Big(\varrho(X^1_s)+\varrho(X^2_s)\Big)\dif s\leq c_3\mE\!\int_0^{\eta}|Z_s|\dif A_1(s).
\end{align*}
Combining the above estimates, and setting
$$
A(t):=A_1(t)+A_2(t),
$$
we arrive at that, for any stopping time $\eta$, there exists a constant $C_0$ such that
\begin{align}
\mE\left[\sup_{t\in[0,\eta]}|Z_{t}|\right]\leq c_0\mE\!\int_0^{\eta}|Z_s|\dif A(s)+c_0\mE\Bigg(\int_0^{\eta}|Z_s|^2\dif A(s)\Bigg)^{\frac{1}{2}}. \label{ee}
\end{align}
By our assumption that $\varrho\in B^0_{q,\infty}(\mR^d)$ with $q>d/\alpha$ and Krylov's estimate (\ref{kry}), we find that for some $c_1>0$,
$$
\mE A_1(t)\leq t+c_1\|\varrho\|_{B^0_{q,\infty}}<\infty.
$$
Since $p>2d/\alpha$, using Fubini's theorem, Krylov's estimate, Minkowski's inequality and taking into account of (\ref{mf}), we can get that for some $c_2>0$,
\begin{align*}
\mE A_2(t)&=\int_{|z|\leq 1}\!\mE\!\int_0^t\Big(\cM|\nabla\cJ_zu|(s,X^1_s)+\cM|\nabla\cJ_zu|(s,X^2_s)\Big)^2\dif s\nu(\dif z)\\
&\leq c_2\!\int_{|z|\leq 1}\!\|(\cM|\nabla\cJ_zu|)^2\|_{L^\infty_TB^0_{p,\infty}}\nu(\dif z)\\
&\leq c_2\!\int_{|z|\leq 1}\!\|(\cM|\nabla\cJ_zu|)^2\|_{L^\infty_p(T)}\nu(\dif z)\leq c_2\!\int_{|z|\leq 1}\!\|\nabla\cJ_zu\|_{L^\infty_p(T)}^2\nu(\dif z),
\end{align*}
where in the second inequality, we have used Krylov's estimate (\ref{kry}) and $p>2d/\alpha$, and in the third inequality we used the fact that
$L^p(\mR^d)\subseteq B^0_{p,\infty}(\mR^d)$.
By our assumption $\beta>1-\alpha/2$, we can always find an $\eps>0$ small enough so that $\beta-\eps>1-\alpha/2$, which in turn means that
$$
2(\alpha+\beta-\eps-1)>\alpha.
$$
Consequently, it follows from \cite[Lemma 2.3]{Zh00} that for some $c_3>0$,
\begin{align*}
\mE A_2(t)\leq c_3\|u\|_{L^\infty_TH^{\alpha+\beta-\eps,p}}^2\int_{|z|\leq 1}|z|^{2(\alpha+\beta-\eps-1)}\nu(\dif z)\leq c_3\|u\|_{L^\infty_TB^{\alpha+\beta}_{p,\infty}}^2<\infty,
\end{align*}
where we have also used the fact that
$B^{\alpha+\beta}_{p,\infty}(\mR^d)\subseteq H^{\alpha+\beta-\eps,p}(\mR^d)$.
Therefore, $t\mapsto A(t)$ is a continuous strictly increasing process. Define for $t\geq 0$ the stopping time
$$
\eta_t:=\inf\{s\geq 0: A(s)\geq t\}.
$$
Then, it is clear that $\eta_t$ is the inverse of $t\mapsto A(t)$. Since $A(t)\geq t$, we further have $\eta_t\leq t$. Plugging $\eta_t$ into (\ref{ee}), we have by a change of variables that
\begin{align*}
&\mE\left[\sup_{s\in[0,t]}|Z_{t\wedge\eta_s}|\right]=\mE\left[\sup_{s\in[0,\eta_t]}|Z_{s}|\right]\\
&\leq c_0\mE\!\int_0^{\eta_t}|Z_s|\dif A(s)+c_0\mE\Bigg(\int_0^{\eta_t}|Z_s|^2\dif A(s)\Bigg)^{\frac{1}{2}}\\
&\leq c_0\mE\!\int_0^{\eta_t}|Z_{s}|\dif A(s)+c_0\mE\Bigg(\int_0^{\eta_t}|Z_{s}|^2\dif A(s)\Bigg)^{\frac{1}{2}}\\
&=c_0\mE\!\int_0^{t}|Z_{t\wedge\eta_s}|\dif s+c_0\mE\Bigg(\int_0^{t}|Z_{t\wedge\eta_s}|^2\dif s\Bigg)^{\frac{1}{2}}\leq c_0\big(t+\sqrt{t}\big)\mE\left[\sup_{s\in[0,t]}|Z_{t\wedge\eta_s}|\right].
\end{align*}
Now taking $t_0$ small enough such that
$$
c_0\big(t_0+\sqrt{t_0}\big)<1,
$$
we get that for almost all $\omega$,
$$
\sup_{s\in[0,\eta_{t_{0}}]}|Z_{s}|=\sup_{s\in[0,t_0]}|Z_{t\wedge\eta_s}|=0.
$$
In particular,
$$
Z_{\eta_{t_0}}=0,\quad a.s..
$$
Repeating the above argument, we can get that for any $k>0$,
$$
\sup_{s\in[0,\eta_{kt_{0}}]}|Z_{s}|=0.
$$
Noticing that $\eta_t$ is strictly increasing, we can get that for all $t\geq 0$,
$$
Z_{t}=0,\quad a.s..
$$
The proof is complete.
\end{proof}

\end{document}